\documentclass{article}

\usepackage[colorlinks=true, allcolors=blue]{hyperref}
\usepackage{amssymb, amsmath, color, xspace, enumerate, mathrsfs,algpseudocode, booktabs, adjustbox, soul, threeparttable, afterpage,longtable, comment, multirow, subcaption, enumerate, longtable, lscape,tabularx, cellspace, svg, tikz, float, wrapstuff}
\usepackage[ruled,vlined,linesnumbered]{algorithm2e}
\usepackage[noblocks,auth-sc]{authblk}
\usepackage[nameinlink]{cleveref}
\usepackage[inline]{enumitem}

    \usepackage[sort,numbers]{natbib}
    \usepackage[margin=1.2in]{geometry}
    \usepackage[justification=centering]{caption}
    \usepackage{amsfonts,amsthm,amsmath}
    \usepackage[title]{appendix}
    \usepackage{charter}
    \newtheorem{theorem}{Theorem}
    \newtheorem{proposition}{Proposition}
    \newtheorem{definition}{Definition}
    \newtheorem{remark}{Remark}
    \newtheorem{example}{Example}
    
\definecolor{lime}{HTML}{A6CE39}
\DeclareRobustCommand{\orcidicon}{
	\begin{tikzpicture} \draw[lime, fill=lime] (0,0) circle [radius=0.16] node[white] { {\fontfamily{qag}\selectfont \tiny ID} };
	\draw[white, fill=white] (-0.0625,0.095) circle [radius=0.007];
	\end{tikzpicture} \hspace{-2mm}
}

\foreach \x in {A, ..., Z}{\expandafter\xdef\csname orcid\x\endcsname{\noexpand\href{https://orcid.org/\csname orcidauthor\x\endcsname}{\noexpand\orcidicon} } }

\setlist[enumerate]{label=(\roman*.)}
\renewcommand{\bar}{\overline}
\renewcommand{\tilde}{\widetilde}
\renewcommand{\hat}{\widehat}

\crefname{lemma}{Lemma}{Lemmata}
\crefname{theorem}{Theorem}{Theorems}
\crefname{claim}{Claim}{Claims}
\crefname{proposition}{Proposition}{Propositions}
\crefname{algorithm}{Algorithm}{Algorithms}
\crefname{equation}{}{}
\crefname{definition}{Definition}{Definition}
\crefname{Cla}{Claim}{Claim}
\crefname{corollary}{Corollary}{Corollaries}
\crefname{remark}{Remark}{Remarks}
\crefname{example}{Example}{Examples}
\crefname{figure}{Figure}{Figures}
\crefname{section}{Section}{Sections}
\crefname{table}{Table}{Tables}
\crefname{enumi}{Statement}{Statements}
\crefname{line}{Step}{Steps}

\usepackage[xindy, acronyms, nowarn]{glossaries}  
\glsdisablehyper
\makeglossaries
\newglossaryentry{IPG}
{
  name={IPG},
  description={Integer Programming Game},
  first={Integer Programming Game (\glsentrytext{IPG})},
  plural={IPGs},
  descriptionplural={Integer Programming Game},
  firstplural={Integer Programming Games (\glsentryplural{IPG})}
}
\newglossaryentry{MNE}
{
  name={MNE},
  description={Mixed Nash Equilibrium},
  first={Mixed Nash Equilibrium (\glsentrytext{MNE})},
  plural={MNEs},
  descriptionplural={Mixed Nash Equilibria},
  firstplural={Mixed Nash Equilibria (\glsentryplural{MNE})}
}
\newglossaryentry{PNE}
{
  name={PNE},
  description={Pure Nash Equilibrium},
  first={Pure Nash Equilibrium (\glsentrytext{PNE})},
  plural={PNEs},
  descriptionplural={Pure Nash Equilibria},
  firstplural={Pure Nash Equilibria (\glsentryplural{PNE})}
}
\newglossaryentry{RBG}
{
  name={RBG},
  description={Reciprocally-Bilinear Game},
  first={Reciprocally-Bilinear Game (\glsentrytext{RBG})},
  plural={RBGs},
  descriptionplural={Reciprocally-Bilinear Games},
  firstplural={Reciprocally-Bilinear Games (\glsentryplural{RBG})}
}
\newglossaryentry{SG}
{
  name={SPG},
  description={Separable-Payoff Game},
  first={Separable-Payoff Game (\glsentrytext{SPG})},
  plural={SPGs},
  descriptionplural={Separable-Payoff Games},
  firstplural={Separable-Payoff Games (\glsentryplural{SG})}
}
\newglossaryentry{PAG}
{
  name={PAG},
  description={Polyhedrally-Approximated Game},
  first={Polyhedrally-Approximated Game (\glsentrytext{PAG})},
  plural={PAGs},
  descriptionplural={Polyhedrally-Approximated Games},
  firstplural={Polyhedrally-Approximated Games (\glsentryplural{PAG})}
}
\newglossaryentry{NASP}
{
  name={NASP},
  description={Nash game among Stackelberg Player},
  first={Nash game among Stackelberg Player (\glsentrytext{NASP})},
  plural={NASPs},
  descriptionplural={Nash games among Stackelberg Players},
  firstplural={Nash games among Stackelberg Players (\glsentryplural{NASP})}
}
\newglossaryentry{NCP}
{
  name={NCP},
  description={Nonlinear Complementarity Problem},
  first={Nonlinear Complementarity Problem (\glsentrytext{NCP})},
  plural={NCPs},
  descriptionplural={Nonlinear Complementarity Problems},
  firstplural={Nonlinear Complementarity Problems (\glsentryplural{NCP})}
}
\newglossaryentry{LCP}
{
  name={LCP},
  description={Linear Complementarity Problem},
  first={Linear Complementarity Problem (\glsentrytext{LCP})},
  plural={LCPs},
  descriptionplural={Linear Complementarity Problems},
  firstplural={Linear Complementarity Problems (\glsentryplural{LCP})}
}
\newacronym{ESGM}{ESGM}{Exhaustive Sample Generation Method}
\newacronym{SGM}{SGM}{Sample Generation Method}
\newacronym{PRLP}{PRLP}{Point-Ray Linear Problem}
\newacronym{CNP}{CnP}{Cut-and-Play}
\newacronym{EO}{ESO}{Enhanced Separation Oracle}
\newacronym{MIP}{MIP}{Mixed-Integer Optimization}

\newcommand{\eg}{\textit{e.g.}}
\newcommand{\ie}{\textit{i.e.}}

\newcommand{\SigmaTwoP}{$\Sigma^p_2$-hard\xspace}
\newcommand{\SigmaTwoPC}{$\Sigma^p_2$-complete\xspace}
\newcommand{\xminusi}{x^{-i}}
\newcommand{\xofi}{x^i}
\newcommand{\tildesigma}{\tilde{\sigma}}
\newcommand{\sigmaofi}{\sigma^i}
\newcommand{\tildesigmaofi}{\tilde{\sigma}^i}
\newcommand{\sigmaminusi}{\sigma^{-i}}
\newcommand{\tildesigmaminusi}{\tilde{\sigma}^{-i}}
\newcommand{\XSetofi}{\mathcal{X}^i}
\newcommand{\XTildeSetofi}{\tilde{\mathcal{X}}^i}
\newcommand{\clconvi}{\cl \conv (\XSetofi)}
\newcommand{\convi}{\conv (\XSetofi)}
\newcommand{\deltai}{\Delta^i}

\newcommand{\R}{\mathbb{R}}

\newcommand{\Z}{\mathbb{Z}}

\newcommand{\ext}{\operatorname{ext}}
\newcommand{\supp}{\operatorname{supp}}
\newcommand{\conv}{\operatorname{conv}}
\newcommand{\cone}{\operatorname{cone}}
\newcommand{\rec}{\operatorname{rec}}
\newcommand{\cl}{\operatorname{cl}}

\newcommand{\theTitle}{The Cut-and-Play Algorithm: Computing Nash Equilibria via Outer Approximations}

\newcommand{\theAbstract}{
	We introduce \emph{Cut-and-Play}, a practically-efficient algorithm for computing Nash equilibria in simultaneous non-cooperative games where players decide via nonconvex and possibly unbounded optimization problems with separable payoff functions. Our algorithm exploits an intrinsic relationship between the equilibria of the original nonconvex game and the ones of a convexified counterpart. In practice, \emph{Cut-and-Play} formulates a series of convex approximations of the game and iteratively refines them with cutting planes and branching operations. Our algorithm does not require convexity or continuity of the player's optimization problems and can be integrated with existing optimization software. We test \emph{Cut-and-Play} on two families of challenging nonconvex games involving discrete decisions and bilevel problems, and we empirically demonstrate that it efficiently computes equilibria while outperforming existing game-specific algorithms.
}

\begin{document}

	\title{\theTitle}
	\date{\vspace{-2em} Margarida Carvalho \orcidA{} Gabriele Dragotto \orcidB{} Andrea Lodi \orcidC{}\\ Sriram Sankaranarayanan  \orcidD{}}
	\author{}

	\maketitle
	\begin{abstract}
		\theAbstract
	\end{abstract}

\section{Introduction}
Decision-making is hardly an individual task; instead, it often involves the mutual interaction of several self-driven decision-markers, or \emph{players}, and their heterogeneous preferences. The most natural framework to model each player's decision problem is often an optimization problem whose solutions, or \emph{strategies}, provide prescriptive recommendations on the best course of action. However, in contrast to optimization, the solution to the game involves the concept of \emph{stability}, \ie, a condition ensuring that the players are playing \emph{mutually-optimal} strategies.
In two seminal papers, \citet{nash_equilibrium_1950,nash_noncoop_1951} formalized a solution concept for finite non-cooperative games, namely, the concept of \emph{Nash equilibrium}.
Nash equilibria are \emph{stable} solutions, as no rational and self-driven decision-maker can \emph{unilaterally} and \emph{profitably} defect them.

Historically, convexity played a central role in shedding light on the existence and computation of Nash equilibria, \eg, see \citet{facchinei_finite-dimensional_2003,daskalakis_non-concave_2022} and the references therein.
Indeed, \citet{vonneumann_zur_1928} proved that 2-player zero-sum games always admit a Nash equilibrium if the players' cost (payoff) functions are convex (concave) in their strategies, and the conditions yielding Nash equilibria are essentially equivalent to linear-programming duality \citep{dantzig1951proof}.
However, the plausibility of the  Nash equilibrium also stems from the availability of algorithms to compute it.  As \citet{roth_economist_2002} argued in his claim “economists as engineers”, computing Nash equilibria plays a central role in designing markets and deriving pragmatic insights. For instance, from the perspective of an external regulator, Nash equilibria are invaluable tools for understanding, designing and intervening in markets where the agents' interests conflict with broader societal objectives.
However, markets and their mathematical models rarely satisfy well-structured convexity assumptions: nonconvexities often model complex operational requirements and are vital to representing reality and extracting valid conclusions. Although there is a wealth of methodologies to compute Nash equilibria for finite and convex games, \eg, normal-form games, little is known about the computation of equilibria in nonconvex settings \citep{daskalakis_non-concave_2022}. This gap represents the core motivation of this work.

This paper presents a practically-efficient method for computing equilibria in non-cooperative games where players solve nonconvex optimization problems.
The nonconvexities may stem from integer variables modeling indivisible quantities and logical conditions, bilevel constraints rendering hierarchical relationships among decision-makers, or from physical phenomena, for instance, water distribution and signal processing \citep{pang_nonconvex_2011,fuller_nonconvex_2022}. Because of their extreme practical interest, several methodologies in optimization have made certain nonconvexities tractable, at least from a computational perspective.
Although some recent papers focused on games with specific nonconvexities, \eg, integer variables as in \citet{koppe_rational_2011,sagratella_computing_2016,cronert_equilibrium_2021,schwarze_branch-and-prune_2022, fuller_nonconvex_2022, carvalho_2020_computing,Dragotto_2021_ZERORegrets,IPGs_2023_Tutorial} and bilevel constraints as in \citet{pang_quasi-variational_2005,WNMS1,sherali1984,Hu2007}, to date, there is no general-purpose algorithm to compute equilibria in games where different nonconvexities arise.

\paragraph{Our Contributions.} This paper presents a practically-efficient algorithm to compute equilibria for a large class of games where players decide by solving optimization problems with nonconvex feasible regions. Specifically, we assume each player optimizes a separable function in its variables, \ie, a function expressed as a weighted sum-of-products function in the player’s variables \citep{dresher_polynomial_1952,dresher_solutions_1953,stein_separable_2008}. We summarize our contributions as follows:

\begin{enumerate}
    \item We prove that when the player's objective functions (\ie, payoffs) are separable, the game admits an equivalent convex representation, even if the players' feasible sets are inherently nonconvex. We present this equivalence in terms of a novel correspondence between the Nash equilibria of the original game and those of a convexified game where each player’s feasible set gets replaced by its convex hull. Furthermore, unlike similar results in optimization, we highlight a fundamental difference between formulating the game over the convex hulls and the \emph{closed} convex hulls of the players’ feasible sets. First, if an equilibrium does not exist in the game formulated over the \emph{closed} convex hulls, it also does not exist in the original game. Second, the existence of an equilibrium in the game formulated over the \emph{closed} convex hulls may not necessarily lead to a corresponding equilibrium in the original game.
    \item We introduce \gls{CNP}, a cutting plane algorithm to compute Nash equilibria for a large family of \emph{simultaneous} and \emph{non-cooperative} $n$-player nonconvex games. In essence, \gls{CNP} exploits a sequence of polyhedral convexifications of the original nonconvex game and involves a well-crafted mix of discrete, \eg, \gls{MIP}, and continuous optimization
    techniques. For instance, it blends concepts such as relaxation (approximation), cutting planes, branching, and complementarity problems. The algorithm finds an exact Nash equilibrium (up to a machine epsilon) or proves its non-existence. To our knowledge, this is also the first approach exploiting outer approximations to compute Nash equilibria.
    \item \gls{CNP} overcomes some well-known algorithmic limitations.
    Specifically, our algorithm \emph{does not}:
          \begin{enumerate*}[label=(\alph*.)]
              \item require that the players' optimization problems are convex or continuous
              \item compute only pure (\ie, deterministic) equilibria
              \item rely on alternative or weaker concepts of equilibria
          \end{enumerate*}. Finally, our algorithm does not require the players to have a bounded set of strategies. On the contrary, it supports unbounded strategy sets and can also certify the non-existence of equilibria. Our primary assumption is that the convex hull of each player's feasible set and its payoff are polyhedral and separable, respectively.
    \item We present an extensive set of computational results on two important families of challenging nonconvex games: \glspl{IPG} and \glspl{NASP}, \ie, a class of simultaneous games among players solving mixed-integer and bilevel optimization problems, respectively. In both cases, our algorithm outperforms the baselines in terms of computing times and social welfare (\ie, the sum of the players' payoffs).
\end{enumerate}

\paragraph{Outline. }  We structure the paper as follows. \cref{sec:background} provides a literature review, and \cref{sec:preliminaries} formalizes the background definitions and our notation. \cref{sec:algo,sec:E0} present the principles behind our approach and introduce the \gls{CNP} algorithm and its practical implementation. \cref{sec:applications} details how to customize \gls{CNP} when some of the game's structure is known, showcasing a comprehensive set of computational results. Finally, we present our conclusions in \cref{sec:conclusions}.

\section{Literature Review}
\label{sec:background}
\citet{nash_equilibrium_1950,nash_noncoop_1951} formalized a concept of stability for non-cooperative simultaneous games through the Nash equilibrium.  We distinguish between two types of equilibria: \glspl{PNE}, where players employ deterministic strategies, and \glspl{MNE}, where players randomize over their strategies. If the game is finite, \ie, there are finitely many players and strategies, Nash proved that an \gls{MNE} always exists. \citet{glicksberg_further_1952} extended the result, proving that an \gls{MNE} always exists when players have continuous payoff functions and compact strategy sets.
Although \glspl{MNE} are guaranteed to exist under some assumptions, an equilibrium may not exist in the general case. Moreover, in \glspl{IPG} and other nonconvex games, deciding if an \gls{MNE} exists is \SigmaTwoPC \citep{vaz_existence_2018,carvalho_2020_computing,WNMS1}, \ie, it is at the second level of the polynomial hierarchy of complexity \citep{woeginger2021}. Besides existence, computing equilibria or certifying their non-existence pose significant algorithmic challenges \citep{daskalakis_complexity_2009}.

\paragraph{Finite Games. }  Most algorithmic approaches for computing Nash equilibria deal with finite games represented in normal form, \ie, through payoff matrices describing the outcome under any combination of the players' strategies.
Besides duality for 2-player zero-sum games \citep{vonneumann_zur_1928}, the first algorithm for computing equilibria in $2$-player normal-form games is the Lemke-Howson algorithm \citep{lemke_equilibrium_1964}.
Although \citet{wilson_computing_1971} and \citet{rosenmuller_generalization_1971} extended the Lemke-Howson algorithm to $n$-player normal-form games, their methods often require the solution of a series of nonlinear systems. More recently, \citet{Sandholm1999} and \citet{porter_simple_2008} proposed two algorithms for 2-player normal-form games exploiting the idea of support enumeration, \ie, the idea of computing an equilibrium by guessing the strategies played with strictly-positive probabilities in an equilibrium. \citet{porter_simple_2008} solve a system of inequalities (nonlinear for more than $2$ players) to determine if a given support (\ie, a subset of pure strategies for each player) leads to an equilibrium. \citet{Sandholm1999} avoid explicit support enumeration by modeling the same idea via a \gls{MIP} problem. In contrast to the above works, we focus on a broader class of games (containing 2-player normal-form games) where the number of pure strategies can be exponential, perhaps uncountable, in the input size of the game.

\paragraph{Continuous Games. } If each player's optimization problem is continuous and convex, equilibrium programming methods can often determine a Nash equilibrium by:
\begin{enumerate*}
    \item reformulating the game as a complementarity or variational inequality problem \citep{cottle_linear_2009}, and
    \item employing globally-convergent Jacobi or Gauss-Seidel algorithms \citep{facchinei_finite-dimensional_2003}.
\end{enumerate*}
On the one hand, these reformulations require restrictive convexity assumptions on the players' optimization problems and may not otherwise guarantee convergence. Besides a few exceptions (\eg, \citet{sagratella_computing_2016}), to date, the majority of equilibrium programming methods require convexity and, otherwise, develop weaker concepts of equilibrium, for instance, \emph{quasi Nash equilibria} \citep{pang_nonconvex_2011,nowak_nonconvex_2021}.
On the other hand, they generally are extremely scalable and efficient under convexity. \gls{CNP} exploits this efficiency by solving, at each step, a convexified game via a complementarity problem.

\paragraph{Integer Nonconvexities. }
In the particular case of integer nonconvexities, \citet{koppe_rational_2011} introduced the taxonomy of \glspl{IPG}, \ie, non-cooperative games where each player solves a parametrized mixed-integer optimization problem. \glspl{IPG} generalize any finite game and implicitly describe the set of strategies via constraints, as opposed to the explicit description of, for instance, normal-form games. Arguably, \glspl{IPG} represent the most prominent emerging family of nonconvex games and have several application domains, for instance, cyber-security \citep{Dragotto_critical_2023}, kidney exchange markets \citep{carvalho_nash_2017}, transportation \citep{sagratella_noncooperative_2020}, and facility location \citep{carvalho2018,LAMAS2018864,cronert_2021,cronert_equilibrium_2021}. Several authors recently developed algorithms to compute \glspl{IPG}' Nash equilibria with a variety of different assumptions \citep{sagratella_computing_2016,schwarze_branch-and-prune_2022,carvalho_2020_computing,cronert_equilibrium_2021,Dragotto_2021_ZERORegrets,harks_2021_generalized}. We refer the reader to the tutorial \citet{IPGs_2023_Tutorial} and the references therein for a detailed survey. Compared to \glspl{IPG} algorithms, \gls{CNP} handles other types of nonconvexities. Furthermore, as opposed to the majority of the algorithms besides the \gls{SGM} of \cite{carvalho_2020_computing} for \glspl{IPG}, \gls{CNP} considers the more general concept of \gls{MNE}, instead of \gls{PNE}.

\section{The Problem and Our Assumptions}
\label{sec:preliminaries}

\subsection{Separable-Payoff Games}
As a standard game-theory notation, let the operator $(\cdot{})^{-i}$ represent the elements of $(\cdot{})$ but the $i$-th element.
We focus on nonconvex \glspl{SG}, a large family of games where each player's objective takes a sum-of-products form as in \cref{def:SG}.

\begin{definition}[Separable-Payoff Game]
    An \gls{SG} $G=(P^1,\dots,P^n)$ is a non-cooperative, complete-information, and simultaneous game among $n$ players where each player $i$ solves
    \begin{align}
        \min_{\xofi}  \{  f^i(\xofi;\xminusi):=(c^i)^\top \xofi + \sum_{j=1}^{m_i} g^i_j(\xminusi) \xofi_j \quad
        \text{ subject to }  \quad\xofi \in \XSetofi \subseteq \R^{m_i}  \}.
        \tag{$P^i$}
        \label{eq:SG}
    \end{align}
    For each player $i$, $\XSetofi \neq \emptyset$ is the set of player $i$'s strategies, $c^i$ is a real-valued vector, and $g^i_j(\xminusi)$ is of the form $g^i_j(\xminusi) = \prod_{\substack{k=1,k\neq i}}^{n} h^i_{j k}(x^k)$ with $h^i_{j k}$ being an affine function.
		An \gls{SG} is \emph{polyhedrally representable} if, for each $i$, $\convi$ (\ie, the convex hull of $\XSetofi$) is a polyhedron.
    \label{def:SG}
\end{definition}
From \cref{def:SG}, the optimization problem $P^i$ of player $i$ is \emph{parametrized} in its opponent choices $\xminusi$.
For each player $i$, we call $\xofi\in\XSetofi$ a \emph{pure strategy}, $\XSetofi$ the feasible set (or set of strategies), and $f^i(\xofi;\xminusi)$ the \emph{payoff} of $i$ under $x=(\xofi,\xminusi)$.\footnote{We slightly abuse the definition of payoff as the players are solving minimization problems. }

\begin{remark}[Linear Form]\label{rem:linearForm}
    Without loss of generality, in \cref{def:SG}, we present \glspl{SG} in a {\em linear (objective) form}, \ie, we let $c^i$ be a vector and $h^i_{jk}(x^k)$ be affine functions.
    If, for any player $i$, $f^i(x^i,x^{-i})=s_i(x^i) + \sum_{j=1}^{m_i} g^i_j(\xminusi) q^i_j(x^i)$, where $s^i$ and $q^i_j$-s are nonlinear functions, and $g^i_j(\xminusi)$ is the product of nonlinear functions $h^i_{jk}(x^k)$ for $k=1,2,\dots,(i-1),(i+1),\dots,n$, we can always reformulate the game so that each payoff has the linear form of \cref{def:SG}.
    To this end, we can introduce
    \begin{enumerate*}
        \item the auxiliary variables $\gamma^i$, $\psi^i_j$ and the constraints $\gamma^i=s^i(\xofi)$, $\psi^i_j = q^i_j(x^i)$ in $P^i$, and
        \item for each $j=1,\dots,m_i$ and $k$, an auxiliary variable $\delta^i_{jk}$ and a constraint $\delta^i_{jk}=h^i_{j k}(x^k)$ in $P^k$.
    \end{enumerate*}
    Thus, the payoff of player $i$ becomes $ \gamma^i + \sum_{j-1}^{m_i} \left( \prod_{k \neq i} \delta^i_{j k} \right) \psi^i_j$, which is in linear form.
    Furthermore, if $s^i(\xofi)$ is convex, we can write the convex inequality $\gamma^i \ge s^i(\xofi)$ instead of the equality. Finally, we remark that \glspl{SG} can represent broad classes of games, for instance, any normal-form game and separable game (which, in addition to separable payoffs, requires $\XSetofi$ to be a compact
    set).
\end{remark}

\paragraph{Mixed Strategies and Equilibria. } For each player $i$, $\sigma^i$ is a \emph{mixed strategy}, or simply a strategy, if it is a probability distribution over the pure strategies $\XSetofi$. Let $\deltai$ be the space of atomic probability distributions over $\XSetofi$ such that $\sigma^i \in \deltai$.  Let $\supp(\sigma^i):=\{\xofi \in \XSetofi : \sigma^i(\xofi)>0 \}$ be the \emph{support} of the strategy $\sigma^i$, where $\sigma^i(\xofi)$ is the probability of playing $\xofi$ in $\sigma^i$.
If a mixed strategy $\sigma^i$ has singleton support, \ie, $|\supp(\sigma^i)|=1$, then it is also a pure strategy.
We denote by $\sigmaminusi \in \prod_{j=1,i\neq j}^n \Delta^j$ the \emph{other players' strategies}, a probability distribution over the strategies of $i$'s opponents.
If  $\XSetofi$ is a compact set, any mixed-strategy $\sigma^i$ of an \gls{SG} has an equivalent finitely-supported mixed-strategy $\hat \sigma^i$ \citep[Theorem 2.8]{stein_separable_2008}. The latter equivalence means that each player $i$'s expected payoff under $\sigma^i$ equals that under $\hat \sigma^i$. This result extends to \glspl{SG}, with $\XSetofi$ being possibly non-compact, as long as the support of $\sigma^i$ is a compact set. This is why, w.l.o.g., $\deltai$ is an atomic distribution over $\XSetofi$. The expected payoff $\mathbb{E}_{\substack{X \sim \sigma}}\left[f^i(X^i;X^{-i})\right]$ for $i$ under $\sigma=(\sigma^1,\dots,\sigma^n)$ is 
\begin{align}
    \quad f^i(\sigma^i;\sigmaminusi) =                                                                    & \sum_{\xofi \in \supp(\sigma^i)}(c^i)^\top \xofi  \sigma^i(\xofi) + \sum_{x \in \supp(\sigma)} \left(\sigma^i(\xofi)\sum_{j=1}^{m_i} g^i_j(\sigma^{-i}(\xminusi)\xminusi) \xofi_j \right).
    \label{eq:expectedPayoff} 
\end{align}
For simplicity, we refer to \cref{eq:expectedPayoff} as $f^i(\sigma^i;\sigmaminusi)$. A strategy $\bar{\sigma}^i$ is a \emph{best response} for player $i$ given its opponents' strategies $\bar{\sigma}^{-i}$ if $f^i(\bar{\sigma}^i;\bar{\sigma}^{-i}) \le f^i(\hat{\sigma}^i;\bar{\sigma}^{-i})$ for any possible \emph{deviation} $\hat{\sigma}^i \in \deltai$. In practice, we can restrict the search of \emph{deviations} $\hat{\sigma}^i$ to pure strategies. A strategy profile $\bar{\sigma}=(\bar{\sigma}^1,\dots,\bar{\sigma}^n)$ is an \gls{MNE} if, for each player $i$ and strategy $\hat{\sigma}^i \in \deltai$, then $f^i(\bar{\sigma}^i;\bar{\sigma}^{-i}) \le f^i(\hat{\sigma}^i;\bar{\sigma}^{-i})$.

\subsection{Polyhedral Representability}
Our algorithmic framework hinges on the assumption of polyhedral representability.
In other words, for any player $i$, we assume that $\convi$, \ie, the convex hull of the set of feasible strategies $\XSetofi$, is a polyhedron.
Whenever $\convi$ is polyhedral for every player $i$, we prove that our algorithm terminates with either an \gls{MNE} or a proof of its non-existence.
The set $\convi$ is a polyhedron if $\XSetofi$ is, for example, a union of finitely many polytopes, the set of mixed-integer points in a polyhedron, or even the union of finitely many polyhedra sharing their set of recession directions. Among the polyhedrally-representable games, we mention the class of \emph{linear} \glspl{IPG} \citep{koppe_rational_2011}, where each player $i$ solves a mixed-integer linear optimization problem.
Besides the assumption of polyhedral representability, our approach is \emph{general} as it does not leverage any game-specific structure, and it computes an \gls{MNE} in any polyhedrally-representable \gls{SG}.

\section{Algorithmic Scheme}\label{sec:algo}

This section outlines \gls{CNP} and its conceptual components. The underlying idea behind our algorithm is to compute an \gls{MNE} by solving a series of convex (outer) approximations of the original game. In principle, solving these convex approximations is computationally more tractable than solving the original nonconvex game. Nevertheless, building an efficient and convergent algorithm poses several algorithmic-design challenges that we discuss in this section.

\paragraph{Notation.} Let $K \subseteq \R^k$ be a closed convex set, and $\rec(K)$ and $\ext(K)$ be the set of its recession directions and extreme points, respectively. An inequality $\pi^\top x \le \pi_0$ is \emph{valid} for $K$ if it holds for any $x \in K$. Given $K$ and a point $\bar x \notin K$, we say $\pi^\top x \le \pi_0$ is a \emph{cut} if it is valid for $K$ and $\pi^\top \bar{x} > \pi_0$. We say $O \subseteq \R^k$ is an outer approximation of $K$ if $K \subseteq O$ and $O$ is polyhedral. Conversely, a (polyhedral) set $I \subseteq \R^k$ is a (polyhedral) inner approximation of $K$ if $I \subseteq K$.

\subsection{Convex Reformulation}
If the players' objectives are separable, we prove the game admits an equivalent convex representation. Specifically, in \cref{thm:convexification}, we establish a correspondence between the \glspl{MNE} of an \gls{SG} instance $G$ and the \glspl{PNE} of a \emph{convexified} instance $\bar G$ where the feasible set for each player's optimization problem is $\convi$ instead of $\XSetofi$. Our result generalizes \citet[Theorem 4]{WNMS1} by letting players have separable payoff functions and arbitrary nonconvex strategy sets. In what follows, we assume that each player's payoff $f^i$ is in the linear form mentioned in \cref{rem:linearForm}.

\begin{theorem} \label{thm:convexification}
    Let $G$ be an \gls{SG} where each player $i$ solves
    $\min_{\xofi}\{f^i(\xofi;\xminusi): \xofi \in \XSetofi\}$.
    Let $\bar G$ be a \emph{convexified} version of $G$ where each player $i$ solves
    $\min_{\xofi}\{f^i(\xofi;\xminusi): \xofi \in \convi \}$.
    For any \gls{PNE} $\bar x$ of $\bar G$, there exists an \gls{MNE} $\hat \sigma$ of $G$ such that, for any player $i$, $f^i(\bar{x}^i; \bar{x}^{-i})=f^i(\hat{\sigma}^i; \hat{\sigma}^{-i})$. Conversely, if $\bar G$ has no \gls{PNE}, then $G$ has no \gls{MNE}.
\end{theorem}

    \begin{proof}[Proof of \cref{thm:convexification}.]
        \label{proof:convexification}
        First, we show that if  $G$ has an \gls{MNE} $\hat \sigma=(\hat \sigma^1,\dots,\hat \sigma^n)$, then the convexified game $\bar G$  has a \gls{PNE} where each player $i$ plays a strategy in $\convi$.
        For each player $i$, we interpret $\hat \sigma^i$ as the mixed (equilibrium) strategy of playing $x^{i^{\ell}}$ with probability $p_\ell$ for $\ell=1,\dots, \kappa_i$.
        We claim that the strategy $\bar{x}^i := \sum_{\ell=1}^{\kappa_i} p_\ell x^{i^{\ell}}$ is feasible in $\bar G$ and a \gls{PNE} of $\bar G$.
        Feasibility follows from the fact that $\bar{x}^i$ is a convex combination of points in $\XSetofi$, and hence $\bar{x}^i \in \convi$.
        Furthermore, the corresponding strategy profile $\bar{x}=(\bar{x}^1,\dots,\bar{x}^n)$ is a \gls{PNE} of $\bar G$ because, for each player $i$, $\bar{x}$ and $\hat{\sigma}$ induce the same payoffs in $\bar G$ and $G$, respectively.
        This is because
        \begin{subequations}
            \begin{align}
                \mathbb{E} _{X \sim \hat \sigma} \left(f^i(X^i;X^{-i})
                \right)
                \quad & =\quad
                \mathbb{E} _{X \sim \hat \sigma} \left(
                (c^i)^\top X^i
                + \sum_{j=1}^{\kappa_i} \left( \prod_{k\neq i} h^i_{j k}(X^k) \right) X^i_j
                \right)
                \label{eq:convexification:1}
                \\
                \quad & =\quad
                \mathbb{E}_{X^i \sim\hat{\sigma}^i} \left(
                (c^i)^\top X^i  \right)
                +\mathbb{E}_{X \sim\hat{\sigma}} \left( \sum_{j=1}^{\kappa_i} \left( \prod_{k\neq i} h^i_{j k}(X^k) \right) X^i_j
                \right)
                \label{eq:convexification:2}
                \\
                \quad & =\quad
                (c^i)^\top\mathbb{E}_{X^i \sim\hat{\sigma}^i} \left(
                X^i  \right)
                +  \sum_{j=1}^{\kappa_i} \left( \prod_{k\neq i}\mathbb{E}_{X^k \sim\hat{\sigma}^k}[ h^i_{j k}(X^k)] \right) \mathbb{E}_{X^i \sim\hat{\sigma}^i} (X^i_j)
                \label{eq:convexification:3}
                \\
                \quad & =\quad
                (c^i)^\top\mathbb{E}_{X^i \sim\hat{\sigma}^i} \left(
                X^i  \right)
                +  \sum_{j=1}^{\kappa_i} \left( \prod_{k\neq i} h^i_{j k}\left (\mathbb{E}_{X^k \sim\hat{\sigma}^k}(X^k)\right ) \right) \mathbb{E}_{X^i \sim\hat{\sigma}^i} (X^i_j)
                \label{eq:convexification:4}
                \\
                \quad & =\quad
                {c^i}^\top \bar{x}^i
                +  \sum_{j=1}^{\kappa_i} \left( \prod_{k\neq i} h^i_{j k}\left (\bar x^k\right ) \right)  \bar x^i_j.
                \label{eq:convexification:5}
            \end{align} \label{eq:convexification}
        \end{subequations}
        Equation \cref{eq:convexification:1} holds because the \gls{SG} is in linear form (see \cref{rem:linearForm}), and \cref{eq:convexification:2} holds because of the linearity of the expectations. Equation \cref{eq:convexification:3} holds because $\hat{\sigma}^{i}$ and $\hat{\sigma}^{k}$ are independent probability distributions of two distinct players $i$ and $k$, and the expectation of the product is the product of the expectations of the random variables.
        Equation \cref{eq:convexification:4} holds due to the linearity of $h^i_{j k}$, and Equation~\cref{eq:convexification:5} holds by the definition of $\bar{x}^i$.
        Thus, for any \gls{MNE} of $G$, a \gls{PNE} of $\bar G$ exists.

        Second, we show that any \gls{PNE} $\bar x$ of $\bar G$ induces an \gls{MNE} $\hat \sigma$ for $G$.
        Because, for each player $i$, $\bar x^i \in \convi$, we rewrite $\bar x^i$ as  $\sum_{\ell=1}^{\kappa_i} p_\ell x^{i^{\ell}}$ for some $p_\ell \ge 0$, $\sum_{\ell=1}^{\kappa_i} p_\ell = 1$ and $x^{i^{\ell}} \in \XSetofi$ for $\ell=1,\dots, \kappa_i$.
        We construct a mixed strategy $\hat \sigma^i$ by letting $i$ select the strategy $x^{i^{\ell}}$ with probability $p_\ell$ for $\ell=1,\dots, \kappa_i$.
        Then, the payoff of $i$ under $\bar x$ is the expression in \cref{eq:convexification:5}.
        By following the equalities \cref{eq:convexification} in the reverse direction, the payoff of $i$ under $\bar x$ is the same as the payoff of $i$ under the mixed strategy $\hat \sigma^i$.
        We claim that $\hat \sigma$ is an \gls{MNE} for $G$ because, for any player $i$ and
        for each unilateral profitable deviation from $\hat \sigma$ in $G$, there exists a unilateral profitable deviation from $\bar x$ in $\bar G$ for $i$.
        Without loss of generality, let the pure strategy $\tilde \sigma^i$ be such deviation for $i$ given $\hat \sigma$ in $G$. Because $\tilde \sigma^i \in \convi$ and \cref{eq:convexification}, $\tilde \sigma^i$ is a deviation for player $i$ given $\bar x$ in $\bar G$, contradicting the assumption that $\bar x$ is a \gls{PNE} for $\bar G$.
        Therefore, any pure strategy in $\bar G$ induces a mixed strategy in $G$ with the same payoff. Then, if $\bar G$ has no \gls{PNE}, then $G$ has no \gls{MNE} (and \emph{vice versa}).
        \end{proof}

\subsection{Building Convex Approximations}
Intuitively, \cref{thm:convexification} proves that any \gls{SG} has an equivalent convex representation where each player $i$ optimizes $f^i(x^i;x^{-i})$ over $\convi$.
Therefore, simultaneously satisfying the optimality conditions of each player's optimization problem in $\bar G$ yields an \gls{MNE} for $G$ and a \gls{PNE} for $\bar G$. Finding a \gls{PNE} in $\bar G$ (if any) is then equivalent to solving an \gls{NCP} expressing the optimality conditions in the form of complementarity conditions.

\paragraph{Approximate Game.}
From a practical perspective, however, the description of each $\convi$ may be challenging to characterize explicitly. For instance, the description of $\convi$ may often be exponentially large in the number of variables or constraints (\eg, mixed-integer sets). Therefore, formulating the game $\bar G$ by employing $\convi$ for each $i$ is practically prohibitive.
Starting from this observation, we devise the concept of \emph{approximate game}, a more tractable convex game where each player $i$'s feasible set is an outer approximation of $\convi$.

\begin{definition}[Approximate Game]
    Let $G$ be an \gls{SG} where each player $i$ solves
    $\min_{\xofi}\{f^i(\xofi;\xminusi): \xofi \in \XSetofi\}$. Then, $\tilde G$ is an \emph{approximate} \gls{SG} of $G$ if each player $i$ solves $\min_{\xofi}\{f^i(\xofi;\xminusi): \xofi \in \XTildeSetofi\}$ with $\XTildeSetofi \supseteq \XSetofi$. Furthermore, $\tilde{G}$ is a \gls{PAG} of $G$ if $\XTildeSetofi$ is a polyhedron for each player $i$.
    \label{def:ApproximateGame}
\end{definition}

\subsubsection{Optimization, Relaxations, and Games.}
In optimization, a feasible solution to the original problem is feasible for its relaxations. However, this relationship may not hold when dealing with games and Nash equilibria. For instance, a game's approximation may admit an \gls{MNE}, whereas the original game may not, or \emph{vice versa}. We illustrate this phenomenon in \cref{ex:RelaxNasp}.

\begin{example}\label{ex:RelaxNasp}
    Consider an \gls{SG} $G$ with $n=2$, where player $1$ solves $\min_x\{\xi x : x\in\R, x\geq 1\}$ and player $2$ solves $\min_\xi\{x\xi  : \xi\in\R, \xi \in [1,2]\}$. This game admits the \gls{MNE} (which is also a \gls{PNE}) $(x, \xi) = (1, 1)$.
    Let $\tilde G$ be \gls{PAG} where the players' feasible regions are $\mathcal{\tilde{X}}^1=\conv(\mathcal{X}^1)$ and $\mathcal{\tilde{X}}^2=\{\xi\in\R: \xi \in [-1,2]\}$, respectively. Although $\tilde G$ has no \gls{MNE}, $G$ admits an \gls{MNE}. If player $1$'s objective changes to $-x\xi$, then $G$ does not have an \gls{MNE}, whereas $\tilde G$ admits the \gls{MNE} $(x, \xi) =(1,-1)$.
\end{example}

An \gls{MNE} for $G$ may not be an \gls{MNE} for one of its \glspl{PAG} as the latter can introduce a destabilizing strategy for $i$ in the associated approximation $\XTildeSetofi$, \ie, a strategy that does not belong to $\convi$ but prevents the existence of that equilibrium in $\tilde G$. This issue is critical when $\XSetofi$ is unbounded or uncountable, as an \gls{MNE} to the original game may not exist.

\subsubsection{Computing Equilibria for the Approximation} In \cref{def:ApproximateGame}, we let $\XTildeSetofi$ be an outer approximation of $\convi$, namely, $\XTildeSetofi$ \emph{enlarges} the feasible set of player $i$. Suppose the approximate game $\tilde G$ is a \gls{PAG}. We can formulate an \gls{NCP} encompassing the optimality conditions associated with each player's optimization problem in $\tilde G$ to determine a \gls{PNE} for $\tilde G$.
At some step $t$ and for some \gls{PAG} $\tilde G$, let $\XTildeSetofi_t=\{\xofi : \tilde{A}^i_t\xofi\le \tilde{b}^i_t, \xofi \ge 0 \}$ be the increasingly-accurate polyhedral approximation of $\convi$ for player $i$. Let $\sigma^i$ and $\mu^i$ be the primal and dual variables of each player's problem in $\tilde G$, respectively. We can compute a \gls{PNE} of $\tilde G$ by solving the \gls{NCP}
\begin{align}
    0 \le \sigmaofi \perp (c^i + \sum_{j=1}^{m_i} g^i_j(\sigmaminusi) + \tilde{A}^{i \top}_t \mu^i) \ge 0, \quad 0 \le \mu^i \perp (\tilde{b}^i_t-\tilde{A}^{i}\sigmaofi) \ge 0 & \qquad  i=1,2,\dots,n,
    \label{eq:NCP}
\end{align}
where $\clubsuit \perp \spadesuit$ is equivalent to $ \clubsuit ^\top \spadesuit=0$. Any solution $\sigma=(\sigma^1,\dots,\sigma^n)$ of \cref{eq:NCP} is a \gls{PNE} for the \gls{PAG}  $\tilde G$ at step $t$. Although a solution to \cref{eq:NCP} includes both $\sigma$ and $\mu$, we omit $\mu$ as we are interested in $\sigma$. If $\XTildeSetofi=\convi$ for any $i$, $\tilde G$ is the exact convex representation of $G$ (\ie, $\tilde G =\bar G$), and the solutions to \cref{eq:NCP} are all the \glspl{MNE} of $G$.

\begin{remark}
	The \gls{NCP} \cref{eq:NCP} carries most of the computational burden of \gls{CNP}, as we empirically show in \cref{sec:applications}. If $g^i_j(\sigmaminusi)$ is linear in $\sigmaminusi$ for any $i$ and $j$, the \gls{NCP} becomes a \gls{LCP}. 
	Unless $P=NP$, there are no polynomial-time algorithms to solve LCPs. However, we can equivalently express an LCP as a MIP and employ efficient MIP solvers. 
	Solving the general \gls{NCP}, however,  often requires mild assumptions on the structure of the underlying game to ensure the convergence of the desired \gls{NCP}'s algorithm  \citep{facchinei_finite-dimensional_2003}.
\end{remark}

\subsection{The Cut-and-Play Algorithm }
\label{sec:sub:CNP}
The workhorse of \gls{CNP} is the \gls{NCP} problem \cref{eq:NCP} associated with each \gls{PAG}. Starting from a game $G$, \gls{CNP} computes the \glspl{PNE} for a finite sequence $t=0,1,2,\dots$ of \glspl{PAG} by
repeatedly solving  \cref{eq:NCP}, and determines whether their solutions are \glspl{MNE} for the original game $G$. If not, the algorithm refines $\XTildeSetofi_t$ for some player $i$ via \emph{cutting} or \emph{branching} on general disjunctions.
The algorithm evokes the same scheme one would use to solve a \gls{MIP} via a branch-and-cut algorithm \citep{padberg_branch-and-cut_1991}, where, instead of a game and an approximate game, one considers a \gls{MIP} and its relaxation.

We emphasize that the approximation of the feasible sets $\XTildeSetofi_t$ is refined via branching and cutting. 
Branching refers to rewriting the feasible set $\XTildeSetofi_t$ as the convex hull of the union of two disjoint sets, such that the union contains every point in $\XSetofi$.  Cutting refers to adding a valid inequality to $\XTildeSetofi_t$. 
Both branching and cutting provide a new outer approximation $\XTildeSetofi_{t+1}$ that is closed and convex.

However, in some cases, $\convi$ may not be a closed set, while its topological closure $\clconvi$ is a polyhedron. 
In this case, the best possible refinement is to get the outer approximation to converge to $\clconvi$ but no further. 
In optimization, the issue of non-closedness is minor, as we can often recover $\varepsilon$-optimal (\ie, approximate) solutions by optimizing over the feasible set's topological closure.
A natural question is if we can use a similar argument for \glspl{SG} and approximate equilibria, \ie, equilibria where each player's strategy is an $\varepsilon$-optimal best response.
In \cref{ex:closure}, we show that, in contrast to optimization, we cannot find an approximate equilibrium for the original game for \emph{any} approximation constant if $\convi$ is not closed.
\begin{example}\label{ex:closure}
    Consider a 2-player \gls{SG} where the players solve 
	\begin{subequations}
		\begin{alignat*}{30}
			\textbf{Player 1:} \qquad & \quad
		\max_{x \in \R^2 } \quad:\quad x_2 \qquad %
		\text{s.t.} %
		 \quad
	x \quad\in\quad \left\{ (x_1, x_2) : x_2 = 0\right\} \cup \left\{ (0,1) \right\} 
	\\
	\textbf{Player 2:} \qquad &
	\max_{y \in \R^2} \quad:\quad  (1-x_1)y_1 + (1-x_2)y_2 
		\end{alignat*}
	\end{subequations}
	Although $\mathcal{X}^1$ is a union of two polyhedra, $\conv(\mathcal{X}^1):=\left\{ (x_1, x_2) \in \R^2 : 0 \le x_2 < 1 \right\} \cup \left\{ (0,1) \right\}$ is not closed. All points along the line $x_2 = 1$ are accumulation points  of $\conv(\mathcal{X}^1)$, but only the point $(0, 1)$ belongs to $\conv(\mathcal{X}^1)$.
	Consider the convexified version of the game. 
	If player $1$ chooses the (infeasible) strategy $(1, 1)$, then, every point in $\R^2$ is feasible and optimal for player $2$ as its objective is $0$. 
	For any other choice of player $1$, player $2$'s objective can be arbitrarily large for an appropriate choice of $y$. 
	However, $(1, 1)$ is not a feasible point  for player $1$. 
	So, for any feasible strategy by player $1$, player $2$ can make its objective arbitrarily large. 
	Thus, the convexified  game has no \gls{PNE} and \cref{thm:convexification} implies that the original game has no \gls{MNE}.
	For a similar reasoning, the game does not even admit an $\varepsilon$-approximate \gls{PNE} for any $\varepsilon > 0$.
	However, $(1,1)$ is a point in the closure of the convex hull of player $1$'s feasible set.
	Thus, if our convexification procedure involves cutting and branching, $(1, 1)$ is always included in the approximate game, and the \gls{PNE} where player 1 plays $(1,1)$ is the output equilibrium. 
	Unlike optimization, there is no sequence of feasible strategies for player $1$, which converges to $(1, 1)$, and is part of any reasonable definition of approximate (pure or mixed-strategy) Nash equilibrium.
\end{example}

A natural question is whether is it possible that the game where each player $i$ plays over $\clconvi$ does not admit an equilibrium while the convex game $\bar G$, \ie, where each player $i$ plays over $\convi$, admits an equilibrium. In \cref{thm:noeq}, we prove this is not possible.
\begin{theorem}
    \label{thm:noeq}
    Let $\bar G$ be a \emph{convexified} \gls{SG} where each player $i$ solves $\min_{\xofi}\{f^i(\xofi;\xminusi): \xofi \in \convi \}$ given some set $\XSetofi$.  Let $\hat G$ be the \gls{SG} where the feasible set of each player in $\bar G$ is replaced with its topological closure, \ie, player $i$ solves
    $\min_{\xofi}\{f^i(\xofi;\xminusi): \xofi \in \clconvi \}$. If $\hat G$ has no \gls{PNE}, then $\bar G$ has no \gls{PNE}.    
\end{theorem}
\begin{proof}
We prove the contrapositive of the statement, \ie, if $\bar G$ has a \gls{PNE}, then $\hat G$ has a \gls{PNE}. In particular, we claim that if $(\bar x^1, \dots, \bar x^n)$ is a \gls{PNE} of $\bar G$, then it is also a \gls{PNE} of $\hat G$. We show this by contradiction.
Suppose $(\bar x^1, \dots, \bar x^n)$ is a \gls{PNE} of $\bar G$, and there exists a player $ i$ and a feasible profitable deviation $\hat x^i \in \clconvi$ such that $f^i(\bar x^i, \bar x^{-i}) - f^i(\hat x^i, \bar x^{-i}) = \epsilon > 0$.
Because $\hat x^i \in \clconvi$, there exists a sequence of points in $\convi$ that converges to $\hat x^i$. By definition, every payoff function $f^i$ is in linear form and thus continuous. Therefore, there exists a point $\tilde x^i \in \convi$ sufficiently close to $\hat x^i$ such that 
$|f^i(\tilde x^i, \bar x^{-i}) - f^i(\hat x^i, \bar x^{-i})| <  \epsilon/2$, which is equivalent to $f^i(\tilde x^i, \bar x^{-i}) - f^i(\hat x^i, \bar x^{-i}) < \epsilon/2$ and $f^i(\hat x^i, \bar x^{-i}) - f^i(\tilde x^i, \bar x^{-i}) < \epsilon/2$. Combining these two inequalities with the definition of $\epsilon$ implies that
$f^i(\bar x^i, \bar x^{-i}) - f^i(\tilde x^i, \bar x^{-i}) \ge \epsilon/2$. However, $\tilde x^i$ is now a feasible profitable deviation in $\bar G$, providing the contradiction that $(\bar x^1, \dots, \bar x^n)$ is not a \gls{PNE} of $\bar G$, and thus completing the proof.
\end{proof}

Given the currently available computational tools, we can only optimize, and thus solve games, over $\clconvi$, and not over $\convi$. 
In \cref{Alg:CnP}, we present a procedure that precisely works with $\clconvi$.
\cref{ex:closure} shows that, unlike optimization problems, working with the closure could lead to important issues. 
We will discuss these issues and possible solutions in \cref{sub:convergence}.

\subsubsection{The Algorithm}
\cref{Alg:CnP} presents the general version of \gls{CNP} for polyhedrally-representable \gls{SG}.
The input of \cref{Alg:CnP} is a polyhedrally-representable \gls{SG} $G$ whereas the output is either an \gls{MNE} $\hat{\sigma}$ or a certificate of its non-existence. 
We will employ an \gls{EO} associated with $P^i$, for each $i$, to separate infeasible equilibria and refine $\tilde G$ via cutting planes, as we will detail in \cref{sec:E0}.
We assume to have access to an initial \gls{PAG} $\tilde{G}$, where, for each player $i$, $\XTildeSetofi_0$ is the starting approximation of the feasible region of player $i$ at step $t=0$. For instance, if $i$ solves a parametrized mixed-integer optimization problem, $\XTildeSetofi_0$ can be its linear relaxation (\ie, the problem without the integrality requirements).
We determine if $\tilde{G}$ has \glspl{PNE} $\tildesigma$ by solving the \gls{NCP} \cref{eq:NCP} induced by $\tilde G$ at step $t$.
\begin{algorithm}[ht]
    \SetKwBlock{Repeat}{repeat}{}
    \DontPrintSemicolon
    \caption{Cut-and-Play for \glspl{SG} \label{Alg:CnP}}
    \KwData{A polyhedrally-representable \gls{SG} instance $G$ }
    \KwResult{
        An \gls{MNE} $\hat{\sigma}=(\hat{\sigma}^1,\dots,\hat{\sigma}^n)$ for $G$ %
        or $\emptyset$ (\ie, no \gls{MNE} exists)}
    $t\gets0,\; \XTildeSetofi_0\gets\{\xofi : \tilde{A}^i\xofi\le \tilde{b}^i, \xofi \ge 0 \}$ for each player $i=1,\ldots,n$\label{Alg:CnP:Relax}\;
    \label{Alg:CnP:While}
    \Repeat{\label{Alg:CnP:Repeat}
        $t\gets t+1$, $\tilde G \gets$ \gls{PAG} where each player $i$ solves $\min_{\xofi} \{ f^i(\xofi;\xminusi) : \xofi \in \XTildeSetofi_t\}$ \label{Alg:CnP:ApproxGame} \;
        $\tilde{\sigma} \gets$ \gls{PNE} deriving from the \gls{NCP} \cref{eq:NCP} \label{Alg:CnP:LCP} \;
        \If{no \gls{PNE} in $\tilde G$\label{Alg:EO:Else}}
        {
            \lIf {$\XTildeSetofi_t = \convi$ for every $i$     \label{Alg:CnP:NoCandidates}}
            {
                \KwRet {$\emptyset$}
            }
            \lElse{
                \textsc{Branch-or-Cut}: refine $\XTildeSetofi_t$ for some $i$ \label{Alg:CnP:Refine}
            }
        }
        \ElseIf{there exists a \gls{PNE} $\tilde{\sigma}$ for $\tilde G$ \label{Alg:CnP:Yes}}
        {
            \For{each player $i=1,2,\dots,n$}{
                $A\gets$\textsc{ESO}$\big(\tildesigmaofi, \XSetofi, f^i((\cdot{}); \tildesigmaminusi)\big)$\label{Alg:CnP:Oracle} \;
                \lIf{$A$ is \texttt{no}}
                {
                    $\XTildeSetofi_{t+1} \gets \XTildeSetofi_t \cap  \{ \xofi : \bar{\pi}^\top \xofi \le \bar{\pi}_0 \} $ \label{Alg:CnP:OracleCut} \tcp*[h]{{\small$\bar \pi,\bar{\pi}_0$ from the \textsc{ESO}}}
                }
            }
            \lIf{$A$ returned \texttt{yes} for every player $i$}
            {
                \textbf{return} $\tildesigma$ \label{Alg:CnP:Return} %
            }
        }
    }
\end{algorithm}

\paragraph{$\tilde G$ has no PNE.} If $\tilde{G}$ has no \gls{PNE}, we cannot infer that $G$ has no \gls{MNE} (see, \eg, \cref{ex:RelaxNasp}), unless $\XTildeSetofi_t = \convi$ for every $i$.
This non-existence can happen when at least one $\XTildeSetofi_t$ is unbounded for some player $i$.
The only viable option is to improve $\tilde{G}$ by refining at least one $\XTildeSetofi_t$ (\cref{Alg:CnP:Refine}). Unfortunately, because $\tilde{G}$ has no \gls{PNE}, we do not have any information on which $\XTildeSetofi_t$ to refine, nor how to refine it. Thus, in \cref{Alg:CnP:Refine} we resort to a \textsc{Branch-or-Cut} subroutine that refines, for some player $i$, $\XTildeSetofi_t$ by either cutting or branching.
If, at step $t$, we need to refine $ \XTildeSetofi_t$ via \cref{Alg:CnP:Refine}, then there exists a $\tildesigmaofi \in \XTildeSetofi_t \backslash \convi$. In a branching refinement, we find two sets $Y^i_{t+1} \subseteq \XTildeSetofi_{t}$ and $Z^i_{t+1}\subseteq \XTildeSetofi_{t}$ such that $\tildesigmaofi \notin \XTildeSetofi_{t+1}=\cl\conv(Y^i_{t+1}\cup Z^i_{t+1})$, with $\XTildeSetofi_{t+1} \subseteq \XTildeSetofi_t$.
This is equivalent to the computation of $\XTildeSetofi_{t+1}$ through Balas' theorem for the union of polyhedra \citep{balas_disjunctive_1998,balas_disjunctive_1985}. Besides \gls{MIP}, branching can handle several other nonconvex problems (\eg, \glspl{LCP}).
In addition (or instead) of branching, \cref{Alg:CnP:Refine} can add to $\XTildeSetofi_{t+1}$ a valid inequality for $\convi$ such that $\XTildeSetofi_{t+1} \subset \XTildeSetofi_{t}$.

\begin{figure}[t]
    \centering
    \includegraphics[width=1\textwidth]{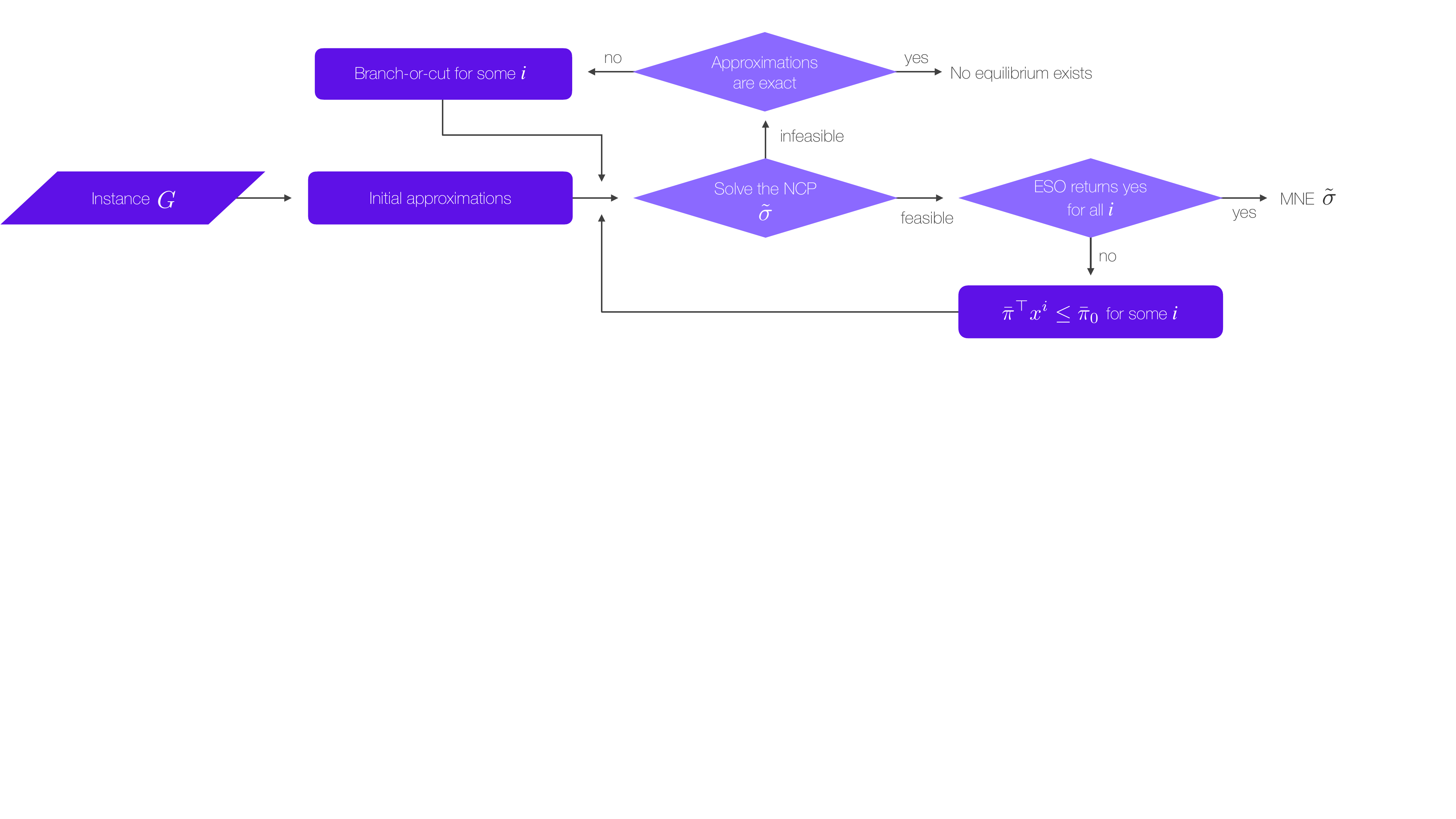}
    \caption{A graphical representation of \gls{CNP}. \label{fig:OA} }
\end{figure}

\paragraph{$\tilde G$ has a PNE.} If $\tilde{G}$ admits a \gls{PNE} $\tilde{\sigma}=(\tilde{\sigma}^1, \dots, \tilde{\sigma}^n)$ (\cref{Alg:CnP:Yes}), then the \gls{EO} will determine whether $\tilde{\sigma}$ is an \gls{MNE} for $G$.
Specifically, given $\tildesigmaofi$ and $\XSetofi$, the \gls{EO} answers $\texttt{yes}$ if $\tildesigmaofi \in \convi$ or $\texttt{no}$ and an hyperplane separating $\tildesigmaofi$ from $\convi$.
On the one hand, if it outputs at least one $\texttt{no}$ for a given player $i$, the oracle certifies that the strategy $\tildesigmaofi$ is infeasible or it is not a best-response to $\tildesigmaminusi$. Then, there exists a valid inequality for $\convi$ that does not hold for $\tildesigmaofi$, \ie, an inequality $\bar{\pi}^\top \xofi \le \bar{\pi}_0$ that refines $\XTildeSetofi_{t}$ to $\XTildeSetofi_{t+i}$.
On the other hand, if the \gls{EO} outputs $\texttt{yes}$ for every player, then $\tildesigma$ is an \gls{MNE} for $G$ (\cref{Alg:CnP:Return}). \cref{fig:OA} illustrates the flow of \gls{CNP}.
\subsubsection{Convergence and Practical Requirements}
\label{sub:convergence}
Keeping \cref{ex:closure} in mind, let us assume that for each player $\convi = \clconvi$; we will discuss the case where this does not hold later. 
To guarantee the convergence of \gls{CNP}, we need to be able to computationally retrieve, in finite time, the description of $\convi$ for any $i$. We formalize this idea with the concept of \emph{computational convexifiability} of \cref{def:convexify}.
\begin{definition}[Computational Convexifiability]
    A set $\mathcal{X}\subseteq \R^k$ is \textit{computationally convexifiable} if 
	\begin{enumerate*} 
		\item $\mathcal{X}$ is closed and convex, or 
		\item the \gls{EO} terminates in a finite number of steps, and given any initial $\tilde {\mathcal{X}} _0\supseteq \mathcal X$, we can obtain $\conv(\mathcal X)$ with a finite sequence of refinements in \cref{Alg:CnP:Refine}. 
	\end{enumerate*}
    \label{def:convexify}
\end{definition}
Naturally, any implementation of an \gls{EO} handling computationally-convexifiable $\convi$ requires the input data of $G$ to be rational. Additionally, branching and cutting should be able to refine the approximations to $\convi$ in a finite number of steps. 
Finally, in \cref{thm:CNP}, we prove that if $\XSetofi$ is computationally convexifiable for every $i$, then \gls{CNP} converges; we provide the full proof in the electronic companion.
\begin{theorem}
	Let $G$ be a polyhedrally-representable \gls{SG}. If $\XSetofi$ is computationally convexifiable for each player $i$, then,  \cref{Alg:CnP} terminates in a finite number of steps and
    \begin{enumerate*}
        \item if it returns $\hat \sigma=(\hat \sigma^1,\dots,\hat \sigma^n)$, then $\hat \sigma$ is an \gls{MNE} for $G$, and
        \item if it returns $\emptyset$, then $G$ has no \gls{MNE}
    \end{enumerate*}.
    \label{thm:CNP}
\end{theorem}

The assumption of computational convexifiability plays a fundamental role in the convergence of CnP (\cref{thm:CNP}). 
Specifically, a set is not computationally convexifiable for two main reasons. 
First, no {\em finite} sequence of branching and cutting can lead to the convex hull, but we could have an infinite sequence of sets, which converges to the convex hull (more formally, the intersection of all these sets is the convex hull).
Second, the convex hull might not be closed as shown in \cref{ex:closure}. 
Hence, no sequence of refinements converging to the convex hull exists.
In either case, \cref{thm:CNP} does not apply. 
The first case is generally due to the fundamental properties of the set $\XSetofi$ itself. 
For example, if $\clconvi$ is not a polyhedron, one might need infinitely many halfspaces to describe it.
This is less of an issue, as refinement can continue indefinitely. 

The second case, however, is more problematic.
Once we obtain $\clconvi$, we cannot refine the outer approximation any further, \ie, to $\convi$.
Nevertheless, we can still verify whether a \gls{PNE} of the outer approximation is feasible, and thus an \gls{MNE} for the original game.
A downside is that there might be infinitely many such equilibria to enumerate. 
While we acknowledge these issues, we also provide reassuring computational evidence.
In the NASPs of \cref{sec:NASP}, the players' feasible sets do not necessarily satisfy $\convi = \clconvi$.
However, we observe that in every instance, the \gls{EO} terminates in a finite number of steps with a \gls{PNE} in the convex hull (and not just in the closure), and the algorithm returns an \gls{MNE} for the original game.
We reiterate that, despite considering $\clconvi$ (as opposed to $\convi$), if we obtain a \gls{PNE} that is in the $\convi$, then we can interpret it as an \gls{MNE} to the original game.
If, however, $\convi \neq \clconvi$ or there is no computational convexifiability, we cannot guarantee that an obtained \gls{PNE} will be in $\convi$.

\section{The Enhanced Separation Oracle}
\label{sec:E0}

Let $\mathcal{X} \subseteq \R^d$ and $\bar{x} \in \R^d$, and assume to have access to an oracle to optimize a linear function over $\mathcal{X}$ in a computationally-tractable manner. The \gls{EO} is an algorithm that, given a point $\bar x$, the set $\mathcal{X}$, and a vector $c\in\R^d$:
\begin{enumerate}
    \item outputs \texttt{yes} and $(V,\alpha)$ if $\bar{x} \in \conv(\mathcal{X})$, with $V \subseteq \mathcal{X}$, and $\alpha \in \R^{|V|}$ being the coefficients of the convex combination of elements in $V$ (\ie, $\bar{x} \in \conv(V) \subseteq \conv (\mathcal{X})$), or
    \item outputs \texttt{no} and a tuple $(\bar{\pi},\bar{\pi}_0)$ so that $\bar{\pi}^\top x \le \bar{\pi}_0 $ is a cutting plane for $\bar{x}$ and $\conv(\mathcal{X})$, \ie, $\bar{\pi}^\top x \le \bar{\pi}_0 $ for any $x \in \conv(\mathcal{X})$ and $\bar{\pi}^\top \bar x > \bar{\pi}_0 $.
\end{enumerate}

Whenever the \gls{EO} outputs \texttt{no}, we separate $\bar x$ from $\conv (\mathcal{X})$ via a cutting plane. Due to \cref{thm:convexification}, this separation task also has the following game-theoretic interpretation if applied to \glspl{SG}; given a set of pure strategies $\XSetofi$ and a point $\tildesigmaofi$, if the \gls{EO} returns \texttt{yes}, then $\tildesigmaofi$ is a mixed strategy, $\supp(\tildesigmaofi)=V$, and $\alpha$ is the vector of probabilities associated with the strategies in $V$.
A theoretical version of this \gls{EO} would include polynomially-many runs of the ellipsoid algorithm, which is theoretically viable yet impractical. 
Furthermore, compared to a standard separation oracle, the \texttt{yes} answer also describes $\bar x$ as a convex combination of points in $\conv(\mathcal{X})$. 
This last requirement is a hard task and further motivates the definition of the \gls{EO}.
Finally, to improve the \gls{EO}'s applicability to \glspl{SG}, we optionally require a vector $c$ to perform an optimization test that provides a sufficient condition for the \gls{EO} to return a $\texttt{no}$ and a \emph{value cut}.

\subsection{Value Cuts}
We start from the concept of \emph{equality of payoffs} \citep{nash_equilibrium_1950,nash_noncoop_1951}, \ie, the concept that, for each player, the payoff of any single pure strategy in the support of an \gls{MNE} strategy must be equal to the \gls{MNE}'s payoff. Formally, let $\sigma^i$ be a (mixed) best-response for player $i$ given $\sigmaminusi$. Then, $f^i(\sigma^i;\sigmaminusi) = f^i(\xofi;\sigmaminusi)$ for any $\xofi \in \supp(\sigma^i)$.
We develop an optimization-based test to diagnose the infeasibility of a given strategy in $\tilde G$ with respect to the original \gls{SG} $G$. Let $\tildesigma=(\tilde \sigma^1,\dots,\tilde \sigma^n)$ for $\tilde G$ be the solution to $\tilde G$  at a given \gls{CNP} step.
Let $\min_{\xofi} \{f^i(\xofi;\tildesigmaminusi): \xofi \in \XSetofi \}$ be the \emph{best-response problem} of $i$ given $\tildesigma^{-i}$, and let $\bar{z}^i$ be its optimal value.
If $\bar{z}^i  > f^i(\tildesigmaofi;\tildesigmaminusi)$, then $\tildesigmaofi \notin \convi$, and a valid separating hyperplane for $\convi$ and $\tildesigma^i$ is $f^i(\xofi; \tildesigmaminusi) \geq \bar{z}^i$. This follows from the equality of payoffs and that $\bar{z}^i$ is the best payoff $i$ can achieve among any pure strategy in $\mathcal{X}^i$. We call these separating hyperplanes \emph{value cuts}. In \cref{thm:ValueCuts}, we prove such inequalities are valid for $\convi$; we provide the proof in the electronic companion.

\begin{proposition}
    Consider an \gls{SG} $G$ and an arbitrary game approximation $\tilde G$ of $G$. Then, for each player $i$ and feasible strategy $\tilde \sigma=(\tildesigma^1,\dots,\tildesigma^n)$ in $\tilde G$,
    $
        f^i(\xofi; \tildesigmaminusi) \geq \inf_{\bar{x}^i} \{f^i(\bar{x}^i;\tildesigmaminusi): \bar{x}^i \in \XSetofi \}
    $
    is a \emph{valid inequality} for $\convi$ if $
        \inf_{\bar{x}^i}\{f^i(\bar{x}^i;\tildesigmaminusi) : \bar{x}^i \in \XSetofi\} = z^i < \infty$.
    If $z^i > f^i(\tildesigma^i;\tildesigma^{-i})$, we call the inequality a \emph{value cut} for $\convi$ and $\tildesigmaofi$.
    \label{thm:ValueCuts}
\end{proposition}

Finally, because we outer approximate $\convi$ for each $i$ and $\tildesigma^i$ generally comes from the \gls{PAG} $\tilde G$, it is not possible to have that $z^i < f^i(\tildesigma^i;\tildesigma^{-i})$. In other words, $\tildesigma^i$ is always a best-response to $\tildesigma^{-i}$ in $\tilde G$, yet, it may be infeasible in $G$.

\subsection{Implementing the Enhanced Separation Oracle}
\label{sec:sub:implementingEO}
We provide an implementation of the \gls{EO} where we require $\conv (\mathcal{X})$ to be a polyhedron. Our implementation decomposes $\conv (\mathcal{X})$ as a conic combination of its rays $\rec(\conv (\mathcal{X}))$ and convex combination of its extreme points $\ext(\conv (\mathcal{X}))$, \ie, it exploits the so-called $\mathcal{V}$-polyhedral representation of $\mathcal{X}$. The \gls{EO} iteratively builds an inner approximation of $\conv (\mathcal{X})$ by keeping track of its rays and vertices. At each \gls{EO}'s call, if the input point $\bar{x}$ cannot be expressed by the incumbent inner approximation of $\conv (\mathcal{X})$, the \gls{EO} either improves the inner approximation by including new vertices and rays or outputs a \texttt{no}. In the case of a \texttt{yes}, this implementation also returns the rays $R$ and the associated conic multipliers $\beta$. In \cref{sub:eliminating}, we also show how to eliminate the conic multipliers and write $\bar x$ exclusively as a convex combination of points $\mathcal{X}$.

\begin{algorithm}[!ht]
    \SetKwBlock{Repeat}{repeat}{}
    \DontPrintSemicolon
    \caption{Enhanced Separation Oracle \label{Alg:EO}}
    \KwData{
        A point $\bar x$, a set $\mathcal{X}$,
	a storage of $V,R$, and optionally a vector $c$}
    \KwResult{ Either: (i.)
        \texttt{yes} and $(V,R,\alpha,\beta)$ if $\bar x \in \conv (\mathcal{X})$, or (ii.) \texttt{no} and a separating hyperplane $\bar{\pi}^\top x \le \bar{\pi}_0$ for $\conv (\mathcal{X}$) and $\bar x$}
    \If{$c$ was provided }
    {
        $\tilde x \gets  \arg \min_{x} \{ c^\top x: x \in \mathcal{X} \}$, with $\bar z = c^\top \tilde x$ \label{Alg:EO:BRValue} \;
        \lIf{$c^\top \bar x < \bar{z} $}
        {
            \KwRet{\texttt{no} and $-c^\top x \le -\bar{z} $} \label{Alg:EO:ValueCut} \hfill %
        }
        \lIf{$\bar x=\tilde x$}
        {
            \KwRet{\texttt{yes} and $(\{\bar{x}\},\emptyset,(1),( ))$}
            \label{Alg:EO:PureBest} %
        }
    }
    \label{Alg:EO:While} \Repeat
    {
        $\mathcal{W} \gets \conv(V) + \cone(R)$. Solve \cref{eq:PRLP} to determine if $\bar x \in \mathcal{W}$ \label{Alg:EO:MinkowskiMembership} \hfill  %

        \lIf{$\bar x \in \mathcal{W}$ }
        {
            \KwRet{\texttt{yes} and $(V,R,\alpha,\beta)$}
            \label{Alg:EO:TrueAnswer}
        }
        \Else(\tcp*[h]{$\bar{\pi}^\top x \leq \bar{\pi}_0 $ is a separating hyperplane for $\bar x$ and $\mathcal{W}$\label{Alg:EO:Farkas}}){
            Let $\mathcal{G}$ be the optimization problem $\max_x \{ \bar{\pi}^\top x : x \in \mathcal{X} \}$ \label{Alg:EO:BlackBox} \;
            \lIf{$\mathcal{G}$ is \emph{unbounded}}
            {
                $R \gets R \cup \{r\}$, where $r$ is an extreme ray of $\mathcal{G}$  \label{Alg:EO:Ray}
            }
            \ElseIf{$\mathcal{G}$ admits an optimal solution $\nu$}{

                \lIf{$\bar{\pi}^\top \nu < \bar{\pi}^\top \bar x$}
                {
                    \KwRet{\texttt{no} and $\bar{\pi}^\top x \le \bar{\pi}^\top \nu  $ \label{Alg:EO:VCut}}
                }
                \lElse{
                    $\nu \gets \arg \max_x \{ \bar{\pi}^\top x : x \in \mathcal{X} \}$, and $V \gets V \cup \{\nu \}$  \label{Alg:EO:Vertex}
                }
            }
        }
    }
\end{algorithm}
\paragraph{The Algorithm. }
\cref{Alg:EO} introduces the implementation of the \gls{EO}. 
This implementation may be warm started with a real-valued vector $c$ to perform the test of \cref{thm:ValueCuts}; for instance, \cref{Alg:CnP:Oracle} of \cref{Alg:CnP} calls the \gls{EO} with $c=c^i+(C^i)^\top \tildesigmaminusi$, the set $\XSetofi$ of $i$, and $\bar x = \tildesigmaofi$. 
As a first step, if $c$ is provided, the algorithm checks if there is any violated value cut by solving the optimization problem $\bar{z}=\min_{x} \{ c^\top x: x \in \mathcal{X} \}$ of \cref{Alg:EO:BRValue}. 
Specifically, the \gls{EO} compares the value of $c^\top \bar x$ to that of $\bar{z}$. 
Let $\tilde x$ be the minimizer yielding $\bar{z}$. 
If  \begin{enumerate*} \item $c^\top \bar x<\bar{z}$, then the \gls{EO} returns a value cut (\cref{Alg:EO:ValueCut}), or \item if the minimizer is $\bar x$ then the \gls{EO} returns \texttt{yes} (\cref{Alg:EO:PureBest}) \end{enumerate*}.
Otherwise, let $V$ and $R$ be a set of vertices and rays of $\mathcal{X}$ that the algorithm can store across its steps. 
We define $\mathcal{W}$ (\cref{Alg:EO:MinkowskiMembership}) as the $\mathcal{V}$-polyhedral inner approximation of $\conv(\mathcal{X})$ such that $\mathcal{W} = \conv(V) + \cone(R)$. 
The central question is then to determine if $\bar x \in \mathcal{W} \subseteq \conv(\mathcal{X})$.

\paragraph{The Point-Ray Separator. }To decide if $\bar x \in \mathcal{W}$, we formulate a linear program expressing $\bar x$ as a convex combination of points in $V$ plus a conic combination of rays in $R$. Let $\alpha$ ($\beta$) be the convex (conic) coefficients for the elements in $V$ ($R$). If $\bar x \in \mathcal{W}$, there exists a solution $(\alpha,\beta)$ to
\begin{align}
    \sum_{v \in V} v \alpha_v + \sum_{r \in R}r \beta_r = \bar x, \quad
    \sum_{v \in V}\alpha_v=1, \quad \alpha \in \mathbb{R}_+^{|V|}, \quad\beta \in \mathbb{R}_+^{|R|}, \label{eq:PrimalPointRay}
\end{align}
where $v$ and $r$ are vertices and rays in $V$, $R$, and $\alpha_v$ and $\beta_r$ are their respective convex and conic multipliers. By duality, \cref{eq:PrimalPointRay} has no solution if there is a $(\bar{\pi},\bar{\pi}_0)$ such that $\bar{\pi}^\top \bar x -\bar{\pi}_0 > 0$, and $v^\top \pi - \pi_0\le 0$, $r^\top \pi \le 0$ for all $v \in V$ and $r \in R$.
Practically, we also maximize the violation $\pi^\top \bar{x} -\pi_0$, and normalize $\pi$ such that $||\pi||_1=1$.  We can equivalently formulate the above requirements as
\begin{subequations}
    \label{eq:PRLP}
    \begin{align}
        \max_{\pi,\pi_0} \qquad  & \pi^\top \bar{x} -\pi_0  \label{eq:PRLP:Obj}    & \qquad\qquad                  & \qquad\qquad\qquad \\
        \text{subject to} \qquad & v^\top \pi  - \pi_0\le 0 \label{eq:PRLP:Vertex} & \forall v \in V, \qquad\qquad & (\alpha)           \\
                                 & r^\top \pi \le 0 \label{eq:PRLP:Ray}            & \forall r \in R, \qquad\qquad & (\beta)            \\
                                 & ||\pi||_1= 1 \label{eq:PRLP:Normaliz}           &                               &
    \end{align}
\end{subequations}
Inspired from \citet{goos_generating_2001} and \citet{chvatal_local_2013}, we define \cref{eq:PRLP} as the \gls{PRLP}. Each vertex $v \in V$ (resp., ray $r \in R$) requires a constraint as in \cref{eq:PRLP:Vertex} (resp., \cref{eq:PRLP:Ray}). As the problem may be unbounded, the normalization constraint \cref{eq:PRLP:Normaliz} truncates the cone of the \gls{PRLP} by requiring the $L_1$-norm of $\pi$ to be $1$.

\begin{figure}[H]
    \captionsetup{}
    \centering
    \begin{subfigure}[b]{0.48\textwidth}
        \centering
        \includegraphics[width=0.7\textwidth]{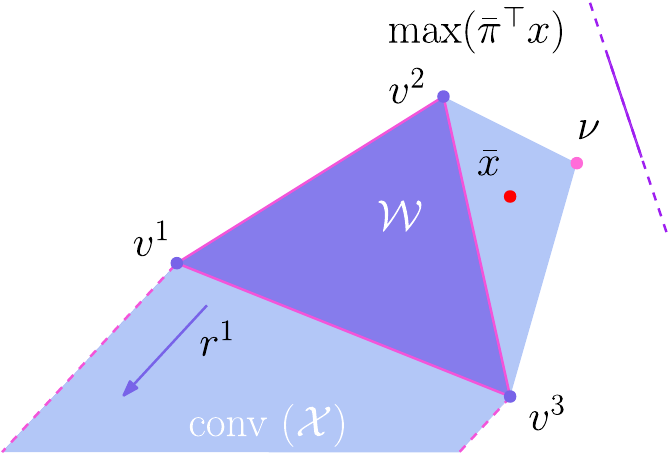}
        \caption{$\nu \notin \mathcal{W}$. Optimizing $\bar{\pi}^\top x$ over $\mathcal{X}$ yields $\nu$. In the next step, the \gls{EO} will return \texttt{yes}.}
        \label{fig:Oracle_Inside}
    \end{subfigure}
    \hfill
    \begin{subfigure}[b]{0.48\textwidth}
        \centering
        \includegraphics[width=0.7\textwidth]{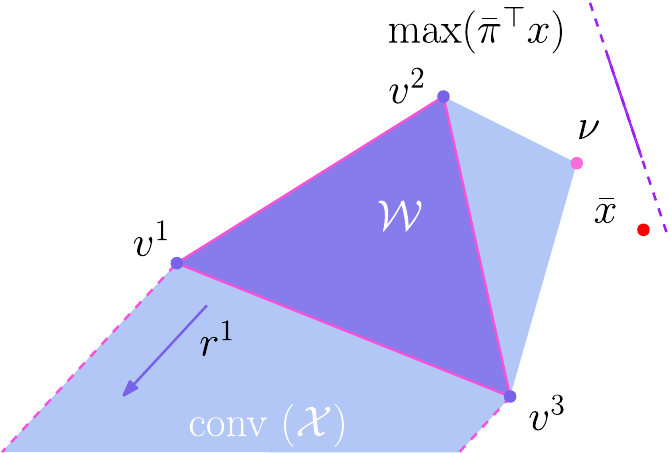}
        \caption{$\nu \notin \mathcal{W}$. Optimizing $\bar{\pi}^\top x$ over $\mathcal{X}$ yields $\nu$. In the next step, the \gls{EO} will return \texttt{no} and $\bar{\pi}^\top \nu < \bar{\pi}^\top \bar x $.}
        \label{fig:Oracle_Outside}
    \end{subfigure}
    \centering
    \caption{A 2-dimensional example of \cref{Alg:EO} separating $\bar x$ from $\conv(\mathcal{X})$. Here, $\mathcal{X}= \{\conv(\{v^2,\nu\})\} \bigcup \{\conv(\{v^1,v^3\})+\cone(\{r^1\})\}$. The set $\conv(\mathcal{X})$ is the light-blue region, whereas  its inner approximation $\mathcal{W}=\conv(\{v^1,v^2,v^3\})$ is in purple.}
    \label{fig:EO_Approx}
\end{figure}

Let $\bar{\pi}, \bar{\pi}_0$ be the optimal values of \cref{eq:PRLP}.
On the one hand, if the optimal objective value of \gls{PRLP} is $0$, the oracle returns \texttt{yes} (\cref{Alg:EO:TrueAnswer}) as $\bar x \in \mathcal{W} \subseteq \conv (\mathcal{X})$.
The convex multipliers $\alpha$ (resp., conic multipliers $\beta$) are the dual values of \cref{eq:PRLP:Vertex} (resp., \cref{eq:PRLP:Ray}).
On the other hand, if $\bar{\pi}^\top \bar{x} - \bar{\pi}_0 >0$, then $\bar{\pi}^\top x \le \bar{\pi}_0$ is a separating hyperplane for $\bar x$ and $\mathcal{W}$. To determine if $\bar{\pi}^\top x \le \bar{\pi}_0 $ is also a separating hyperplane for $\bar x$ and $\conv (\mathcal{X})$, the \gls{EO} solves the problem $\mathcal{G}=\max_x \{ \bar{\pi}^\top x : x \in \mathcal{X}\}$ (\cref{Alg:EO:BlackBox}). If $\mathcal{G}$ is unbounded, then its extreme ray $r$ is a new ray for the set $R$. Conversely, if $\mathcal{G}$ admits an optimal solution $\nu$, the latter is a new vertex for the set $V$ (\cref{Alg:EO:Vertex}). Furthermore, if $\bar{\pi}^\top \nu<\bar{\pi}^\top \bar x $, then $\bar x$ is infeasible. In practice, this means $\bar x$ is separated from $\conv (\mathcal{X})$ by $\bar{\pi}^\top x \le \bar{\pi}^\top \nu $, and the \gls{EO} returns \texttt{no}. If this is not the case, the \gls{EO} identified a new vertex (or ray), and the process restarts from \cref{Alg:EO:MinkowskiMembership}. We represent \cref{Alg:EO} in \cref{fig:EO_Approx}.

\paragraph{Practical Considerations. } First, similarly to \citet{goos_generating_2001}, we can modify \cref{Alg:EO:BlackBox} of \cref{Alg:EO} to retrieve multiple vertices and rays violating $\bar{\pi}^\top x \le \bar{\pi}_0 $, and subsequently add them in \cref{Alg:EO:Vertex} and \cref{Alg:EO:Ray}. In this way, the inner approximation  $\mathcal{W}$ tends to build faster without significantly impacting the computational overhead. Second, the normalizations of the \gls{PRLP} in \cref{eq:PRLP} are practically pivotal as they affect the algorithm's overall stability (and convergence) through the generated cutting planes. Because normalizations tend to significantly affect the generators' performance \citep{dey_approximating_2015,fischetti_separation_2011,bixby_solving_2002,goos_generating_2001}, we normalize \cref{eq:PRLP} with \cref{eq:PRLP:Normaliz}.
Finally, in \cref{prop:EOTerminates}, we show that our implementation of the \gls{EO} terminates in a finite number of steps; we defer the proof to the electronic companion.
\begin{theorem}
    The \gls{EO} terminates in a finite number of steps if $\conv (\mathcal{X})$ is a polyhedron.
    \label{prop:EOTerminates}
\end{theorem}

\subsection{Eliminating the Conic Coefficients}
\label{sub:eliminating}
In \cref{thm:convexification}, we interpret any {\em convex} combination of strategies in $\mathcal{X}^i$ in a game-theoretic fashion, \ie, as a mixed strategy for player $i$, where each element in the convex combination is a pure strategy whose probability of being played is given by the associated convex combination's coefficient.
However, \cref{Alg:EO} returns, whenever the answer is \texttt{yes},
a \emph{proof of inclusion} $(V,R,\alpha,\beta)$ that also includes a conic combination of the extreme rays in $R$.
As the game-theoretic interpretability of the solution given by \gls{CNP} is critical, especially when players solve unbounded problems, we provide a simple algorithm to \emph{repair} a proof of inclusion $(V,R,\alpha,\beta)$ to a proof of inclusion $(V,\alpha)$ that \emph{does not} include any conic combination. We illustrate its intuition in \cref{example:conic}.

\begin{example}[Conic Combinations]
    Consider the example in \cref{fig:Repair}, and assume that \cref{Alg:EO} returns \texttt{yes}, and a proof of inclusion $(V,R,\alpha,\beta)$.
    Let $V$ be $\{ v^1, v^2 \}$, $R$ be $\{ r^1 \}$, and $\mathcal{W}$ be $\conv(V) + \cone(R)$. Let $\mathcal{X}$ be made by the points $\{v^1,v^2,v^3\}$ and the sequence of points $v^{3,1},v^{3,2},\dots$ along $r^1$. Although the proof of inclusion of $\bar x$ employs the ray $r^1$, we can equivalently express $\bar x$ as a convex combination of $v^1$, $v^2$ and $v^{3,1}$ without resorting to $r^1$.
    \begin{wrapstuff}[type=figure,width=.4\textwidth]
        \centering
        \includegraphics[width=\textwidth]{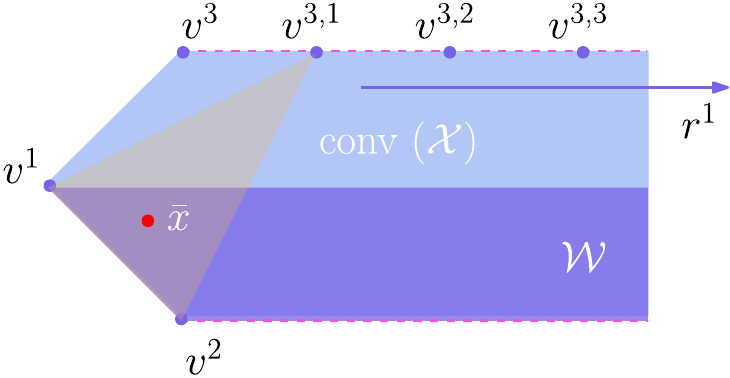}
        \caption{Rewriting $\bar x$ as a convex combination of $v^1$, $v^2$ and $v^{3,2}$ (orange area).}
        \label{fig:Repair}
    \end{wrapstuff}
    \label{example:conic}
\end{example}

\paragraph{A Repair Algorithm.} If \gls{CNP} terminates with an \gls{MNE} that includes rays, \cref{Alg:Repair} provides the repairing routine that eliminates the rays from the proof of inclusion of each player $i$.
\cref{Alg:Repair} requires the point $\bar x \in \conv(\mathcal{X})$, the set $\mathcal{X}$,  an arbitrarily-large constant $B \in \R$ and the proof of inclusion $(V,R,\alpha,\beta)$ from the \gls{EO} of \cref{Alg:EO}.
\cref{Alg:Repair} iteratively attempts to express $\bar x$ as a convex combination of points in $V$ by augmenting $V$ with some points (not necessarily extreme points) of $\mathcal{X}$.
Initialize $\tilde{V}$ as $V$. At each step, the algorithm augments $\tilde{V}$ with the optimal solutions of $\max_x \{ r^\top x : x \in \mathcal{X},  r^\top x \le B \}$ for any $r \in R$ (\cref{Alg:Repair:Vertex}).
If $\bar x \in \tilde{V}$ due to the \gls{PRLP} of \cref{Alg:Repair:MinkowskiMembership}, then the algorithm returns $(\tilde{V},\tilde\alpha)$, where $\tilde \alpha$ are the dual variables of the \gls{PRLP}.
Otherwise, the algorithm increments $B$ with a positive integer and keeps iterating. We remark that \cref{Alg:Repair} terminates in a finite number of steps, as $x \in \conv(\mathcal{X})$, and there exists a $B'$ such that $r^\top x \le B'$ for any $r \in R$ and $x \in \mathcal{X}$.

\begin{algorithm}[!ht]
    \SetKwBlock{Repeat}{repeat}{}
    \DontPrintSemicolon
    \caption{Repairing Algorithm \label{Alg:Repair}}
    \KwData{
        A point $\bar x$, $\mathcal{X}$, %
         a large $B \in \R$, and $(V,R,\alpha,\beta)$ from \cref{Alg:EO}}
    \KwResult{A proof of inclusion $(V,\alpha)$}
    $\tilde{V} \gets V$\;
    \Repeat
    {
        \lFor{$r \in R$}{
            $\nu \gets \arg \max_x \{ r^\top x : x \in \mathcal{X},  r^\top x \le B \}$, and $\tilde{V} \gets \tilde{V} \cup \{\nu \}$  \label{Alg:Repair:Vertex}
        }
        Solve the \gls{PRLP} \cref{eq:PRLP} to determine if $\bar x \in \conv(\tilde{V})$ \label{Alg:Repair:MinkowskiMembership} \hfill \;
        \lIf{$\bar x \in \conv(\tilde{V})$}
        {
            \KwRet{\texttt{yes} and $(\tilde{V},\tilde \alpha)$}
            \label{Alg:Repair:TrueAnswer}
            \textbf{else }increase $B$
        }
    }
\end{algorithm}

\section{Applications}
\label{sec:applications}
In this section, we evaluate \gls{CNP} on two challenging nonconvex games, demonstrate the algorithm's effectiveness, and establish a solid computational benchmark against the literature. We consider \glspl{IPG} and \glspl{NASP}, two games among players solving integer and bilevel problems. In both cases, determining the existence of an equilibrium is generally \SigmaTwoP. Thus, we expect the computation of an equilibrium, or the determination of its non-existence, to be challenging.
We compare \gls{CNP} against the most advanced problem-specific algorithms and empirically demonstrate that our algorithm is scalable, practically efficient, and can exploit problem-specific structures.

\paragraph{Reciprocally-bilinear Games.} In our tests, we focus on the \glspl{RBG}, a subclass of \gls{SG} in the form of \cref{def:RBG}.

\begin{definition}[Reciprocally-Bilinear Game]
    An \gls{RBG} is an \gls{SG} where, for each player $i$, $f^i(\xofi;\xminusi):=(c^i) ^\top \xofi + (\xminusi)^\top C^i\xofi$, where $C^i$ and $c^i$ are a matrix and vector with rational entries.
    \label{def:RBG}
\end{definition}
In \glspl{RBG}, determining if a \gls{PAG} $\tilde G$ admits an \gls{MNE} is equivalent to solving a \gls{LCP}. Both specialized \gls{LCP} solvers (\eg, \emph{PATH} from \citet{dirkse_path_1995,ferris_interfaces_1999}), and \gls{MIP} reformulations \citep{kleinert_why_2020} can solve \glspl{LCP}.
Although a \gls{MIP} reformulation does not exploit the underlying complementarity structure, \gls{MIP} solvers can optimize several computationally-tractable functions over the set of the \gls{LCP}'s solutions and, therefore, select an \gls{MNE} in the \gls{PAG} $\tilde G$ that maximizes a given objective function. In this sense, \gls{CNP} supports heuristic equilibria selection and enables the user to select the preferred balance between equilibrium quality (\eg, given a function that measures its quality) with the time required for its computation.
In our tests, we employ \emph{Gurobi} 9.2 as \gls{MIP} solver and \emph{PATH} as an \gls{LCP} solver. \footnote{Our tests run on an \emph{Intel Xeon Gold 6142} with $128$GB of RAM. The code and the instances are available in the package \emph{ZERO} \citep{Dragotto_2021_ZEROSoftware} at \url{https://ds4dm.github.io/ZERO}.} We report all the time-related results as shifted geometric means with a shift of $10$ seconds.

\subsection{Integer Programming Games}
We consider a class of \glspl{IPG} where each player $i$ solves the mixed-integer optimization problem
\begin{align}
    \max_{\xofi} \big\{ (c^i)^\top \xofi + (\xminusi)^\top C^i \xofi \quad \text{subject to} \quad A^i\xofi \le b^i,  \quad \xofi\ge 0, \quad \xofi_j \in \Z \quad  \forall j \in \mathcal{I}^i \big\}.
    \label{eq:IPG}
\end{align}
\noindent The matrix $A^i$ and the vector $b^i$ have rational entries for any $i$, and $\mathcal{I}^i$ contains the indexes of the integer-constrained variables. Reciprocally-bilinear \glspl{IPG} have applications in several domains, \eg, revenue management, healthcare, cybersecurity, and sustainability \citep{cronert_equilibrium_2021,Dragotto_critical_2023,carvalho_nash_2017}, and can also represent any 2-player normal-form game.
We also note that for each player $\convi = \clconvi$ and the feasible sets are computationally convexifiable. 
Thus, the issues described in \cref{ex:closure} do not apply to \glspl{IPG}.

\subsubsection{Customizing \gls{CNP}. }
We tailor \cref{Alg:CnP} as follows:
\begin{enumerate}
    \item
		A set of inequalities describing $\convi$ is a \emph{perfect formulation} of $\XSetofi$ \citep{conforti_extended_2010}.
		  In \cref{Alg:CnP:Refine} of \cref{Alg:CnP}, we can refine each player's approximation $\XTildeSetofi$ with any family of \gls{MIP} inequalities. In our tests, we employ Gomory Mixed-Integer (\emph{GMIs}), Mixed-Integer Rounding (\emph{MIRs}), and Knapsack Cover (\emph{KPs}) inequalities through the software \emph{CoinOR Cgl} \citep{CoinCGL}.
          Furthermore, at step $t=0$ of \gls{CNP}, we let each $\XTildeSetofi_{0}$ be the associated linear relaxation.
    \item If \cref{Alg:CnP:Oracle} of \cref{Alg:CnP} returns \texttt{no} at step $t$, we add additional cuts (\eg, \emph{GMIs}) to refine $\XTildeSetofi_{t+1}$. Furthermore, whenever the value cuts exhibit a numerically-ill behavior, we attempt to substitute it with a \gls{MIP} cut with better numerical properties, \eg, we generate a \gls{MIP} cut that cuts off the incumbent solution.
\end{enumerate}

\begin{remark}[Branching may not be enough]
    The \gls{EO}'s cuts in \cref{Alg:CnP:OracleCut} of \cref{Alg:CnP} are essential for the convergence of \gls{CNP}. Specifically,
    a pure branching algorithm (\ie, without cutting)
    cannot guarantee the separation of an infeasible strategy $\tildesigmaofi$ of $\tilde G$ from $\convi$.
    Consider, for instance, the incumbent \gls{PNE} $\tildesigma$ of $\tilde G$ at step $t$:
    the associated strategy $\tildesigmaofi$ of player $i$ may not be an extreme point of $\XTildeSetofi_t$ and, on the contrary,
    it can be in the interior of $\XTildeSetofi_t$. In this case, the branching operation may not be able to exclude $\tildesigmaofi$ from the refined $\XTildeSetofi_{t+1}$ (\cref{fig:Branching}).
    Indeed, finding the branching strategy excluding infeasible strategies may be a hard problem. Thus, a pure and naive branching algorithm might not separate the incumbent \gls{PNE} of $\tilde G$ from the players' approximations.
    \begin{figure}[!ht]
        \centering              \includegraphics[width=0.80\textwidth]{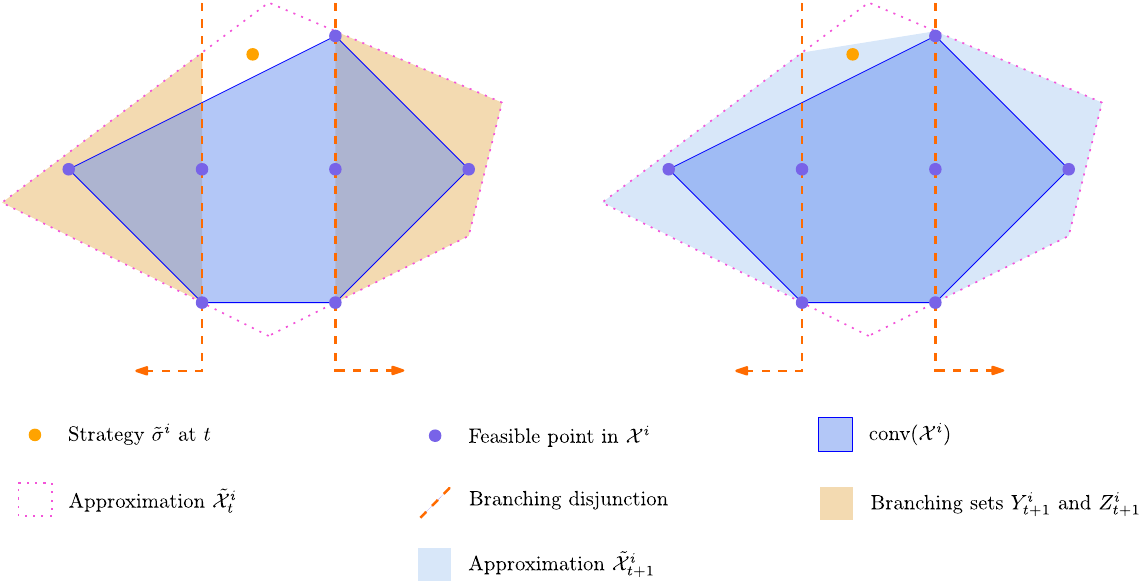}
        \caption{At step $t$, branching might not exclude an infeasible \gls{PNE} from $\XTildeSetofi_{t+1}$. Given $\tildesigmaofi$ in the interior of $\XTildeSetofi_{t}$, the branching operation results in a refined $\XTildeSetofi_{t+1}=\cl \conv(Y^i_{t+1} \cup Z^i_{t+1})$ that does not exclude $\tildesigmaofi$.}
        \label{fig:Branching}
    \end{figure}
    Finally, whereas branching cannot always separate a strategy from the \emph{closure} of the convex hull, neither cutting nor branching can separate a strategy in $\clconvi$  from $\convi$.
\end{remark}

\subsubsection{Computational Tests}
\paragraph{Instances and Parameters.} We compare \gls{CNP} with the \gls{SGM} algorithm of \citet{carvalho_2020_computing}, currently the most efficient algorithm to compute an \gls{MNE} in \glspl{IPG}. We employ the popular instances of the \emph{knapsack game}, where each player $i$ solves a knapsack problem with $m$ items; namely, each player $i$ has a set of strategies $\mathcal{X}^i=\big \{w^i \xofi \le \bar{w}^i, \; \; \xofi \in \{0,1\}^m\big\}$.
The parameters $c^i \in \Z^{m}$ and $w^i \in \Z^{m}_+$ are integer vectors representing the profits and weights of player $i$, respectively; The parameter $\bar{w}^i$ is the knapsack capacity, whereas $C^i \in \Z^{(n-1)m \times m}$ is the concatenation of $(n-1)$ diagonal matrices. The elements on the diagonals are the so-called \emph{interaction coefficients} associated with each of the $(n-1)$ other players in the game and their $m$ decision items (in the lexicographic order given by each player's index).  Thus, players interact only when selecting corresponding items, \ie, $\xofi_j=x^o_j=1$ for $i,o \in \{1,\dots,n\}$. We remark that $C^i$ are different among players, and their entries are integer-valued but not necessarily positive, \ie, the interaction for a given item can be positive or negative.
Because $\convi$ and $\XTildeSetofi$ are compact for every $i$, any \gls{PAG} admits a \gls{PNE} and \gls{CNP} purely acts as a cutting plane algorithm. Furthermore, although an equilibrium always exists, computing it is at least $PPAD$-hard.
We solve the \glspl{LCP} in \cref{eq:NCP} with either: \begin{enumerate*} \item \emph{PATH}, thus computing \emph{a} feasible \gls{MNE} for each \gls{PAG}, or \item \emph{MIP}, by optimizing the quadratic \emph{social welfare} function $SW(\sigma)=\sum_{i=1}^nf^i(\sigma^i;\sigma^{-i})$. \end{enumerate*}
When optimizing the social welfare, we aim to find equilibria exhibiting desirable properties from the perspective of a regulator, \ie, we aim to heuristically select the equilibria that favor the \emph{society} the most. Finally, we set a numerical tolerence of   $\varepsilon = 3\times 10^{-5}$, a time limit of $300$ seconds, and we employ the $70$ instances from \citet{carvalho_2020_computing}, where $n\in \{2,3\}$ and $10\leq m \leq 100$.

\paragraph{Results. } \cref{tab:IPGs} provides an overview of the results by categorizing the instances in two sets: the small instances (\ie, $mn\le80$) in rows $2-8$, and the large ones (\ie, $mn>80$) in rows $9-15$. The first column reports the algorithm, where, for \gls{CNP}, we specify whether we use \emph{Gurobi} or \emph{PATH} (resp., \emph{CnP-MIP} and \emph{CnP-PATH}). Column \emph{O} is the objective type, either \emph{F} for feasibility or \emph{Q} when \gls{CNP} optimizes $SW$ while solving each \gls{PAG}. Column \emph{C} reports the aggressiveness of the additional \gls{MIP} inequalities generated. Specifically \emph{C} can be:
\begin{enumerate*}
    \item $-1$ if \gls{CNP} do not use \gls{MIP} cuts, or
    \item $0$ if it adds \gls{MIP} cuts to replace numerically-ill value cuts, or
    \item $1$ if it concurrently adds \gls{MIP} cuts on top of numerically-stable value cuts.
\end{enumerate*}
The columns \emph{Time} and \emph{\#TL} report the time (in seconds) and the number of time limit hits, respectively. Column \emph{MinMax} reports the absolute difference (in seconds) between the maximum and the minimum computing time. Column \emph{SW\%} reports the average social welfare improvement compared to the \gls{MNE} computed by \gls{SGM}. Finally, we report the average number of steps (\emph{\#It}), cuts added (\emph{Cuts}), and \emph{MIP} cuts (\emph{MIPCuts}). We remark that  \emph{Cuts} includes any cut from the \gls{EO} (including value cuts) and the \gls{MIP} cuts.

\paragraph{Discussion.} \gls{CNP} always computes equilibria with remarkably modest computing times. Furthermore, \gls{CNP} improves the average social welfare compared to \gls{SGM}. This is mainly due to the algorithm's approximation structure: whereas \gls{SGM} approximates the player's strategy sets from the inside, \gls{CNP} outer approximates them and thus has a larger search space. The social welfare improves when \gls{CNP} leverages a \gls{MIP} solver, with the average welfare values almost doubling. However, this improvement significantly increases the computing times as \gls{MIP} solvers cannot exploit the structure of \cref{eq:NCP}. When \gls{CNP} uses \emph{PATH}, there are dramatic computing time improvements over the whole set of instances.
Furthermore, the more \gls{MIP} cuts, the fewer steps are required to converge to an \gls{MNE}. When combining \gls{CNP} with \gls{MIP} cuts, we generally observe a reduction in the computing time and a lower \emph{MinMax} value. Finally, we remark that although our instances are as large as the ones considered in \citet{Dragotto_2021_ZERORegrets,carvalho_2020_computing,cronert_equilibrium_2021}, \gls{CNP} exhibits limited computing times.
\begin{table}[ht]
    \centering
        \resizebox{0.95\textwidth}{!}{
        \begin{tabular}{l@{\hspace{1em}}l@{\hspace{1em}}l@{\hspace{3em}}r@{\hspace{3em}}rr@{\hspace{2em}}r@{\hspace{2em}}r@{\hspace{2em}}r@{\hspace{2em}}r@{\hspace{1.2em}}}
            \toprule
            \textbf{Algorithm} & \textbf{O} & \textbf{C} & \textbf{Time (s)} & \textbf{MinMax} & \textbf{\#TL} & \textbf{SW\%} & \textbf{\#It} & \textbf{Cuts} & \textbf{MIPCuts} \\
            \toprule
            \textbf{SGM}       &            & \textbf{}  & 0.73              & 21.43           & 0             & 0.0\%         & 8.43          & -             & -                \\
            \textbf{CnP-MIP}   & Q          & -1         & 6.58              & 287.52          & 0             & 13.5\%        & 7.80          & 9.57          & 0.00             \\
            \textbf{CnP-MIP}   & Q          & 0          & 6.13              & 287.01          & 0             & 12.9\%        & 5.73          & 6.47          & 2.30             \\
            \textbf{CnP-MIP}   & Q          & 1          & 6.31              & 287.52          & 0             & 13.3\%        & 3.50          & 9.60          & 7.47             \\
            \textbf{CnP-PATH}  & F          & -1         & 0.36              & 10.54           & 0             & 1.8\%         & 7.60          & 10.2          & 0.00             \\
            \textbf{CnP-PATH}  & F          & 0          & 0.05              & 0.19            & 0             & 2.9\%         & 5.27          & 5.90          & 2.07             \\
            \textbf{CnP-PATH}  & F          & 1          & 0.04              & 0.19            & 0             & 4.9\%         & 3.23          & 8.87          & 7.10             \\
            \hline
            \textbf{SGM}       &            & \textbf{}  & 20.86             & 300.00          & 6             & 0.0\%         & 18.58         & -             & -                \\
            \textbf{CnP-MIP}   & Q          & -1         & 61.08             & 294.50          & 0             & 22.5\%        & 13.70         & 17.00         & 0.00             \\
            \textbf{CnP-MIP}   & Q          & 0          & 57.85             & 299.45          & 1             & 19.6\%        & 11.62         & 12.62         & 3.45             \\
            \textbf{CnP-MIP}   & Q          & 1          & 68.20             & 299.04          & 0             & 22.3\%        & 9.48          & 16.80         & 10.32            \\
            \textbf{CnP-PATH}  & F          & -1         & 6.68              & 80.89           & 0             & 15.7\%        & 13.55         & 16.35         & 0.00             \\
            \textbf{CnP-PATH}  & F          & 0          & 4.48              & 74.37           & 0             & 15.7\%        & 9.62          & 10.25         & 2.42             \\
            \textbf{CnP-PATH}  & F          & 1          & 4.32              & 75.88           & 0             & 15.9\%        & 8.22          & 14.35         & 8.43
            \\
            \bottomrule
        \end{tabular}
    }
    \caption{\glspl{IPG} results on the knapsack game. Extended results are available in the electronic companion. \label{tab:IPGs}}
\end{table}
\subsection{Nash Games Among Stackelberg Players}
\label{sec:NASP}
\glspl{NASP} \citep{WNMS1} are \glspl{SG} where each player $i$ solves a bilevel program. 
This family of games is instrumental in energy and pricing, as it represents complex systems of hierarchical interaction and combines simultaneous and sequential interactions. In a \gls{NASP}, each player $i$'s feasible region is
\begin{subequations}
    \begin{align}
        \XSetofi := \left \{ x^i : A^i \xofi \le b^i, \quad \xofi \ge 0, \quad\; \xofi=(w^i,\hat{y}^i), \quad \hat{y}^{i} \in \text{SOL}(w^i) \right \} \label{eq:Stackelberg:Lowerlevel},
    \end{align}
\end{subequations}
where $\xofi$ is partitioned into the leader's variables $w^i$, and the followers' variables $y^i$. For each player $i$, there are $u_i \in \Z$ followers controlling the variables $y^{i,k}$ with $k=1,\dots,u_i$. Each follower $k$ solves a convex quadratic optimization problem in $y^{i,k}$ parametrized in $w^i$. The set $\text{SOL}(w^i)$ in \cref{eq:Stackelberg:Lowerlevel} represents the solutions $y^i=(y^{i,1},\dots,y^{i,u_i})$ to the $i$-th player lower-level simultaneous game.
Therefore, the feasible set $\XSetofi$ for each player is a set of linear constraints \cref{eq:Stackelberg:Lowerlevel} plus the optimality of the followers' game \cref{eq:Stackelberg:Lowerlevel}.
Each player's feasible region $\XSetofi$ is a union of finitely many polyhedra \citep{Basu2018a} and we can express $\XSetofi$ as
\begin{align}
    \underbrace{
        \{
        (\xofi,z^i):A^i\xofi\le b^i, z^i=M^i\xofi+q^i, \xofi \ge 0, z^i \ge 0
        \}
    }_{\XSetofi_0}
    \bigcap_{j \in \mathcal{C}^i}( \{ (\xofi,z^i):z^i_j=0 \} \cup \{(\xofi,z^i):\xofi_j=0 \}), \label{eq:UnionRef}
\end{align}
where $\mathcal{C}^i$ is a set of indexes for the complementarity equations that equivalently express $\text{SOL}(w^i)$. For any player $i$, we refer to $\XSetofi_0$ in \cref{eq:UnionRef} as its \emph{polyhedral relaxation}. In other words, $\XSetofi_0$ is the polyhedron containing the leader constraints, the definitions for $z^i$, and the non-negativity constraints. Because $\convi$ may be unbounded, an \gls{MNE} in \glspl{NASP} might not exist. Indeed, the problem of determining if a \gls{NASP} admits an equilibrium is \SigmaTwoP \citep{WNMS1}. Finally, $\convi$ is not necessarily closed because of the union of polyhedra appearing in \cref{eq:UnionRef}. Therefore, \gls{CNP} actually works with $\clconvi$.

Moreover, the players' feasible sets in \glspl{NASP} may not be computationally convexifiable. 
Since $\convi$ may not be closed, we may encounter some of the issues illustrated in \cref{ex:closure}. 
This is an example of a problem class where all our assumptions for finite termination, as per \cref{thm:CNP}, do not necessarily hold. 
However, we try \cref{Alg:CnP} on this problem, and assess the performance.

\subsubsection{Customizing \gls{CNP}. } We express each player's optimization problem by reformulating $\XSetofi$ as in \cref{eq:UnionRef}.
In \gls{CNP}, the initial relaxation of \cref{Alg:CnP:Relax} is the polyhedral relaxation $\XTildeSetofi_0$ for each player $i$, \ie, the feasible region where we omit \emph{all} the complementarity constraints. As the description of \cref{eq:UnionRef} only needs a finite number of complementarity conditions in $\mathcal{C}^i$, the branching steps account for finding the disjunction $j$ associated with the complementarity condition $z^i_j=0$ or $\xofi_j=0$. Let $J^i_t$ be the set of disjunctions added at step $t$ of the algorithm. Then, \cref{Alg:CnP:Refine} will include in $\XTildeSetofi_{t+1}$ at least one complementarity $j$ for some player $i$ so that $j \in \mathcal{C}^i \backslash J^i_t$. This boils down to the computation of $\XTildeSetofi_{t+1}$ as the union of $\XTildeSetofi_{t} \cap \{ (\xofi,z^i) : \xofi_{j^i}\le 0\}$ and $\XTildeSetofi_{t} \cap \{ (\xofi,z^i) : z^i_{j^i}\le 0\}$, where $\xofi_{j^i}$,$z^i_{j^i}$ are the terms involved in the $j$-th complementarity of $i$. If, at some step $t$, $J^i_t = \mathcal{C}^i$ for all $i$, then $\convi = \XTildeSetofi_{t}$ and the algorithm terminates.
Furthermore, we develop two custom branching rules for \glspl{NASP}. Assume a \gls{PAG} admits an \gls{MNE} $\tildesigma$ at step $t-1$ and needs to determine the best branching candidate in $\mathcal{C}^i \backslash J^i_t$ at step $t$. Then, we branch via:
\begin{enumerate}
    \item \emph{Hybrid branching. } For any candidate set $\mathcal{C}^i \backslash J^i_t$, we select the candidate $j$ that minimizes the distance between $\tildesigmaofi$ and the set $\XTildeSetofi_{t+1}$ that includes the $j$-th complementarity. Therefore, we select $j$ by solving
          $\min_{\lambda^i}
              \lbrace
              (\lambda^i)^2 : \tildesigmaofi \in \cl \conv \left(
              \lbrace
              \XTildeSetofi_{t} \cap
              \lbrace (\xofi,z^i) : \xofi_{j^i}\le 0 \rbrace
              \rbrace
              \cup
              \lbrace
              \XTildeSetofi_{t} \cap
              \lbrace (\xofi,z^i) : z^i_{j^i}\le 0 \rbrace
              \rbrace
              \right)$,
          where $\lambda$ is the vector of violations associated with each constraint in $\XTildeSetofi_{t}$ plus the disjunctions on the $j$-th complementarity.
    \item \emph{Deviation branching.} We solve the best-response problem of $i$ given $\tildesigmaminusi$, and compute an optimal solution $\tilde{x}^i$. We select the first (given an arbitrary order) candidate $j$ that encodes the polyhedron containing $\tilde{x}^i$.
\end{enumerate}

\subsubsection{Computational Tests }
\paragraph{Instance and Parameters.} We set a numerical tolerance of $\varepsilon = 10^{-5}$, and consider a time limit of $300$ seconds. We employ the $50$ instances \emph{InstanceSet B} from \citet{WNMS1}, where each instance has $7$ players with up to $3$ followers each, and compare against the problem-specific (sequential) \emph{inner approximation} algorithm (\emph{Inn-S}) from \citet{WNMS1}. We also introduce $50$ harder instances \emph{H7} with $7$ players with $7$ followers each. Large \glspl{NASP} instances, such as the \emph{H7} set, are numerically badly scaled and thus helpful to perform stress tests on the numerical stability of the algorithms.

\paragraph{Results and Discussion.} 
First, we note that, despite the fact that the convex hulls of the feasible sets may not be closed, in every instance, we obtained a \gls{PNE} to the outer approximation that is in the convex hull of the feasible set, 
\ie, in the computational experiments we never encountered the type of issue described in \cref{ex:closure}.
Thus, for every instance we either identified an \gls{MNE} or proved its non-existence computationally. 
We report the aggregated results in \cref{tab:NASPS1}. 
The first column reports the type of algorithm: the baseline (\emph{Inn-S}) or \gls{CNP} with either the hybrid (\emph{HB}) or deviation (\emph{DB}) branching. 
Furthermore, we test  \emph{Inn-S-1} and \emph{Inn-S-3}, two configurations of the inner approximation algorithm (we refer to the original paper for a detailed description). 
The second column reports the instance set (\emph{Inst}), \ie, either \emph{B} or \emph{H7}. 
The three subsequent pairs of columns report the average computing time (\emph{Time}) and the number of instances (\emph{\#}) for which the algorithm either: \begin{enumerate*} \item found an \gls{MNE} (\emph{EQ}), or \item proved that no \gls{MNE} exists (\emph{NO\_EQ}), or \item terminated with either an \gls{MNE} or a proof of its non-existence without exhibiting numerical issues (\emph{ALL}) \end{enumerate*}. 
The last two columns report the number of numerical issues (\emph{\#NI}) and time limit hits (\emph{\#TL}) each algorithm encountered. 
The baseline \emph{Inn-S-1} systematically fails on the set \emph{H7} due to the size of the descriptions of $\convi$. 
Indeed, \emph{Inn} exhibits significant numerical issues in the set \emph{H7}, even though the former is a problem-specific algorithm. 
On the contrary, \gls{CNP} performs consistently and is especially effective in the hard set \emph{H7}.
The running times of both algorithms are comparable in the instance set \emph{B}, and \gls{CNP} is competitive with \emph{Inn} while not being a problem-specific algorithm.

\begin{table}[!ht]
    \centering
    \caption{\glspl{NASP} results. \label{tab:NASPS1}}
        \resizebox{0.85\textwidth}{!}{
        \begin{tabular}{ll@{\hspace{4em}}r@{\hspace{1em}}r@{\hspace{2em}}r@{\hspace{1em}}r@{\hspace{2em}}r@{\hspace{1em}}r@{\hspace{3em}}r@{\hspace{1em}}r}
            \toprule
            \textbf{Algo}    & \textbf{Inst} & \textbf{Time (s)}      & \textbf{\#}                & \textbf{Time (s)}       & \textbf{\#} & \textbf{Time (s)} & \textbf{\#N} & \textbf{\#NI} & \textbf{\#TL} \\
            \hline
                             &               & \multicolumn{2}{c}{EQ} & \multicolumn{2}{c}{NO\_EQ} & \multicolumn{2}{c}{ALL} &                                                                                \\
            \hline
            \textbf{Inn-S-1} & \textbf{B}    & 6.22                   & 49                         & 69.76                   & 1           & 6.56              & 50           & 0             & 0             \\
            \textbf{Inn-S-3} & \textbf{B}    & 4.94                   & 49                         & 23.96                   & 1           & 5.12              & 50           & 0             & 0             \\
            \textbf{CnP-HB}  & \textbf{B}    & 7.47                   & 46                         & 29.37                   & 1           & 7.71              & 47           & 3             & 0             \\
            \textbf{CnP-DB}  & \textbf{B}    & 9.45                   & 46                         & 11.81                   & 1           & 9.50              & 47           & 3             & 0             \\
            \hline
            \textbf{Inn-S-1} & \textbf{H7}   & -                      & 0                          & -                       & 0           & 300.00            & 46           & 4             & 46            \\
            \textbf{Inn-S-3} & \textbf{H7}   & -                      & 0                          & -                       & 0           & -                 & 0            & 50            & 0             \\
            \textbf{CnP-HB}  & \textbf{H7}   & 53.79                  & 41                         & -                       & 0           & 73.45             & 50           & 0             & 9             \\
            \textbf{CnP-DB}  & \textbf{H7}   & 52.58                  & 35                         & -                       & 0           & 88.92             & 50           & 0             & 15            \\
            \bottomrule
        \end{tabular}
    }
\end{table}

\section{Concluding Remarks}
\label{sec:conclusions}
In this work, we presented \gls{CNP}, a practically-efficient algorithm for computing Nash equilibria in \glspl{SG}, a large class of non-cooperative games where players solve nonconvex optimization problems. We showed that nonconvex \glspl{SG} admit equivalent convex formulations where the players' feasible sets are the convex hulls of their original nonconvex feasible regions. 
Starting from this result, we designed \gls{CNP}, a cutting-plane algorithm to compute \glspl{MNE} for polyhedrally-representable \glspl{SG} or certify their non-existence. We defined the concept of \emph{game approximation}, and we employed it through \gls{CNP} by building an increasingly accurate sequence of convex approximations converging to an equilibrium or a certificate of its non-existence. Our algorithm is general and exhibits solid computational performance. Although \gls{CNP} is a general-purpose algorithm, we also demonstrated how to tailor it to exploit the structure of specific classes of games. In addition, we also provided a series of insights into the relationships between equilibria, approximations, and approximated equilibria.

Given the generality of our algorithm, we prudently believe improvement opportunities lie ahead. We hope our contribution can inspire future methodological development in equilibria computation in nonconvex games.
Among those, we foresee an extension of our methodology beyond the assumption of polyhedral representability.
Finally, as Nash equilibria play a pivotal role in designing and regulating economic markets, we hope our algorithm will enable economists and optimizers to design complex markets stemming from sophisticated economic models.

\section*{Acknowledgement}
We thank Federico Bobbio, Didier Chételat, Aleksandr Kazachkov, Rosario Scatamacchia, and Mathieu Tanneau for the insightful discussions concerning this work. Gabriele Dragotto, Andrea Lodi, and Sriram Sankaranarayanan are thankful for the support of the Canada Excellence Research Chair in “Data Science for Real-time Decision-making” at Polytechnique Montreal. Finally, we are deeply thankful to the two anonymous referees for their constructive comments and suggestions.

\bibliography{Main.bib}
\label{sec:bib}

\appendix
\section{Proof of Theorem \ref{thm:CNP}.}
	\begin{proof}[Proof of \cref{thm:CNP}.]
		\emph{Finite Termination.} We show that \gls{CNP} terminates in a finite number of steps. The calls to solve the \gls{NCP} in \cref{Alg:CnP:LCP}, and to the \gls{EO} in \cref{Alg:CnP:Oracle} terminate in a finite number of steps because we assume the players' feasible sets are computationally-convexifiable. The only loop that could potentially not terminate is the \texttt{repeat} starting in \cref{Alg:CnP:Repeat}.

		First, we restrict to the case where the set $\XSetofi$ is bounded or finite for any player $i$. Any \gls{PAG} $\tilde G$ necessarily admits a \gls{PNE} if the approximations $\XTildeSetofi_t$ are compact. Therefore, the algorithm never triggers \cref{Alg:EO:Else}. Thus, \cref{Alg:CnP:OracleCut} is the only step refining the sets $\XTildeSetofi_t$ for some player $i$. As $\XSetofi$ is computationally convexifiable for every player $i$, the \gls{EO} can, in a finite number of steps, refine the approximations $\XTildeSetofi_t$ to $\convi$. As a consequence, the algorithm converges (in the worst case) to $\bar G$ (\ie, the exact convex approximation), and the correctness of the resulting \gls{MNE} follows from \cref{thm:convexification}.

		Second, if $\XSetofi$ is unbounded for some player $i$, then a $\gls{PNE}$ for a given \gls{PAG} $\tilde G$ may not exist, and the algorithm may enter \cref{Alg:EO:Else}. However, because $\XSetofi$ is computationally convexifiable for every player $i$, the branching step and the \gls{EO} can refine, in a finite number of steps, the approximations $\XTildeSetofi_t$ to $\convi$. Therefore, even in the unbounded case, the algorithm correctly returns an \gls{MNE}, similarly to the bounded case.

		\emph{Proof of statements (i) and (ii).}
		We show that $\hat \sigma$ is an \gls{MNE} for $G$.
		If the algorithm returns $\hat \sigma$, there exists an approximate game $\tilde{G}$ in \cref{Alg:CnP:ApproxGame} with a \gls{PNE} $\tilde \sigma=\hat \sigma$. Let the step associated with $\tilde G$ having a \gls{PNE} $\hat \sigma$ be denoted with $t=\theta$, and, for each player $i$, let $\XTildeSetofi_\theta$ be the associated feasible region $i$. Then, for any player $i$ and $\bar \sigma^i \in \XTildeSetofi_\theta$, it follows that $f^i(\hat {\sigma}^i;\hat{\sigma}^{-i}) \le f^i(\bar \sigma^i;\hat{\sigma}^{-i})$, \ie,
		no player $i$ has the incentive to deviate from $\hat {\sigma}^i$ to any other strategy $\bar \sigma^i \in \XTildeSetofi_\theta$ in $\tilde G$. Because $\convi \subseteq \XTildeSetofi_\theta$ for any $i$, then the previous inequality holds for any $\bar \sigma^i \in \convi$. Moreover, by construction, $\hat {\sigma}^i \in \convi$; otherwise, the \gls{EO} would have returned \texttt{no} for player $i$.
		\end{proof}

\section{Proof of Theorem \ref{prop:EOTerminates}.}
	\begin{proof}[Proof of \cref{prop:EOTerminates}.]
		The \gls{EO} inner approximates the polyhedron $\conv (\mathcal{X})$ with its $\mathcal{V}$-representation, which is made of finitely many extreme rays and vertices. Hence, we have to prove that the \gls{EO} never finds, at any step, a vertex $\nu$ (ray $r$) in \cref{Alg:EO:Vertex} (\cref{Alg:EO:Ray}) so that $\nu$ is already in $V$ ($r$ is already in $R$). This implies that the repeat loop in \cref{Alg:EO} terminates.

		The inequality after the else statement in \cref{Alg:EO:Farkas} is valid for $\mathcal{W}$ if and only if $v^\top \bar{\pi} \le \bar{\pi}_0 $ for any $v \in V$, and $r^\top \bar{\pi} \le 0$ for any $r \in R$ as of \cref{eq:PRLP:Vertex} and \cref{eq:PRLP:Ray}. Also, because the latter inequality is a separating hyperplane between $\mathcal{W}$ and $\bar x$, then $\bar{\pi}^\top \bar x > \bar{\pi}_0 $. However, it may not necessarily be a valid inequality for any element in $\ext(\conv (\mathcal{X}))$ and $\rec(\conv (\mathcal{X}))$. Therefore, we must consider the optimization problem $\mathcal{G}$ in \cref{Alg:EO:BlackBox}. On the one hand, if $\mathcal{G}$ is bounded, let $\nu$ be its optimal solution. Then, either \begin{enumerate*} \item $\bar{\pi}^\top \nu < \bar{\pi}^\top \bar x $, with $\bar{\pi}^\top x \leq \bar{\pi}^\top \nu $ being a separating hyperplane between $\conv (\mathcal{X})$ and $\bar x$, and the algorithm terminates and returns \texttt{no}, or \item  $\bar{\pi}^\top \nu \ge \bar{\pi}^\top \bar x $, $\nu$ is necessarily a vertex of $\ext(\conv (\mathcal{X})) \backslash V$ violating $\bar{\pi}^\top x \le \bar{\pi}_0 $, and the algorithm updates $V \gets V \cup \{\nu\}$ \end{enumerate*}. On the other hand, if $\mathcal{G}$ is unbounded, then there exists an extreme ray $r$ so that $r^\top \bar{\pi}> 0$. Then, $r$ is necessarily in $\rec(\conv (\mathcal{X})) \backslash R$, $\bar{\pi}^\top r > \bar{\pi}_0 $, and the algorithm updates $R \gets R \cup \{r\}$ and returns to \cref{Alg:EO:MinkowskiMembership}. As there are finitely many extreme rays and vertices, the algorithm terminates.
		\end{proof}

\section{Proof of Proposition \ref{thm:ValueCuts}}
	\begin{proof}[Proof of \cref{thm:ValueCuts}.]
		\label{proof:ValueCuts}
		If the infimum is attained at a finite value $z^i$, this implies that player $i$ cannot achieve a payoff strictly less than $z^i$ given the other players' strategies $\tildesigmaminusi$; in other words, $z^i$ is the payoff associated with the best response of $i$ to $\tildesigmaminusi$. Consider the problem
		\begin{equation}
			\inf_{\bar{x}^i}\{f^i(\bar{x}^i;\tildesigmaminusi) : \bar{x}^i \in \convi\}.
			\label{eq:proof:value_min}
		\end{equation}
		The above problem attains a finite infimum because $z^i$ is finite.
		Let $\bar{v}^i$ be an optimal solution, \ie, a best response, and $\bar{z}^i$ be the finite optimal value of \cref{eq:proof:value_min}, respectively. We claim that $z^i=\bar{z}^i$ for any $\bar{v}^i \in \convi$ that solves \cref{eq:proof:value_min}.
		On the one hand, assume that $\bar{z}^i < z^i$. Note that $\bar{v}^i \in \convi$ is a convex combination of points in $\XSetofi$, and $f^i(x^i;\tildesigmaminusi)$ is linear in $x^i$. Then, any point $\hat v^i$ involved in the convex combination resulting in $\bar{v}^i$ belongs to $\XSetofi$; this also imply that $z^i$ is not the optimal value of $\inf_{\bar{x}^i}\{f^i(\bar{x}^i;\tildesigmaminusi) : \bar{x}^i \in \XSetofi\}$, resulting in a contradiction. On the other hand, assume that $\bar{z}^i > z^i$. Because the solutions to $\inf_{\bar{x}^i}\{f^i(\bar{x}^i;\tildesigmaminusi) : \bar{x}^i \in \XSetofi\}$ are also feasible for \cref{eq:proof:value_min}, $\bar v^i$ cannot be a minimizer of \cref{eq:proof:value_min}.
		\end{proof}

\section{IPGs Results}
\cref{tab:IPGRes,tab:IPGRes2} presents the full computational results for our experiments. The column names are analogous to those of \cref{tab:IPGs}, with the addition of a few columns.  Specifically, we report the value of the social welfare in \emph{SW} and the average numbers of: \begin{enumerate*}
	\item cuts from the \gls{EO} excluding value cuts (\emph{ESOCuts}),
	\item value cuts from the \gls{EO} (\emph{VCuts}).
\end{enumerate*}
Finally, in the time column, we report in parenthesis the time the algorithm spent to compute the first \gls{MNE}; this is relevant when \gls{CNP} optimizes the social welfare function via a \gls{MIP} solver.

\newpage

\begin{table}[H]
	\centering
	\caption{\glspl{IPG} complete results, first set.  \label{tab:IPGRes}}
		\resizebox{0.95\textwidth}{!}{
		\begin{tabular}{lll@{\hspace{4em}}rrrrrrrr}
			\toprule
			\textbf{Algorithm} & \textbf{O} & \textbf{C} & \textbf{Time (s)} & \textbf{\#TL} & \textbf{SW} & \textbf{\#It} & \textbf{Cuts} & \textbf{ESOCuts} & \textbf{VCuts} & \textbf{MIPCuts} \\
			\hline
			\texttt{n=3 m=10}  &            &            &                   &               &             &               &               &                  &                &                  \\
			\hline
			\textbf{SGM}       & -          & -          & 2.11              & 0             & 632.99      & 10.00         & -             & -                & -              & -                \\
			\textbf{CnP-MIP}   & Q          & -1         & 0.47 (0.23)       & 0             & 812.48      & 4.50          & 5.0           & 2.0              & 3.0            & 0.0              \\
			\textbf{CnP-MIP}   & Q          & 0          & 0.31 (0.14)       & 0             & 812.98      & 4.60          & 4.8           & 2.0              & 1.1            & 1.7              \\
			\textbf{CnP-MIP}   & Q          & 1          & 0.20 (0.08)       & 0             & 820.71      & 2.60          & 7.2           & 0.5              & 1.1            & 5.6              \\
			\textbf{CnP-PATH}  & F          & -1         & 0.02              & 0             & 706.66      & 5.00          & 5.9           & 2.0              & 3.9            & 0.0              \\
			\textbf{CnP-PATH}  & F          & 0          & 0.02              & 0             & 718.13      & 4.50          & 4.9           & 2.0              & 1.5            & 1.4              \\
			\textbf{CnP-PATH}  & F          & 1          & 0.03              & 0             & 742.87      & 2.00          & 5.4           & 0.3              & 0.7            & 4.4              \\
			\hline
			\texttt{n=2 m=20}  &            &            &                   &               &             &               &               &                  &                &                  \\
			\hline
			\textbf{SGM}       & -          & -          & 0.01              & 0             & 658.31      & 5.40          & -             & -                & -              & -                \\
			\textbf{CnP-MIP}   & Q          & -1         & 0.96 (0.25)       & 0             & 684.19      & 6.40          & 6.3           & 4.4              & 1.9            & 0.0              \\
			\textbf{CnP-MIP}   & Q          & 0          & 0.93 (0.29)       & 0             & 683.91      & 6.10          & 5.9           & 3.0              & 1.2            & 1.7              \\
			\textbf{CnP-MIP}   & Q          & 1          & 0.75 (0.18)       & 0             & 682.69      & 3.70          & 7.6           & 1.4              & 0.9            & 5.3              \\
			\textbf{CnP-PATH}  & F          & -1         & 0.05              & 0             & 645.44      & 5.30          & 5.5           & 3.1              & 2.4            & 0.0              \\
			\textbf{CnP-PATH}  & F          & 0          & 0.04              & 0             & 664.44      & 4.90          & 4.7           & 1.8              & 1.2            & 1.7              \\
			\textbf{CnP-PATH}  & F          & 1          & 0.03              & 0             & 656.44      & 3.10          & 6.2           & 1.2              & 0.4            & 4.6              \\
			\hline
			\texttt{n=3 m=20}  &            &            &                   &               &             &               &               &                  &                &                  \\
			\hline
			\textbf{SGM}       & -          & -          & 0.20              & 0             & 1339.98     & 9.90          & -             & -                & -              & -                \\
			\textbf{CnP-MIP}   & Q          & -1         & 29.74 (1.49)      & 0             & 1488.96     & 12.50         & 17.4          & 7.0              & 10.4           & 0.0              \\
			\textbf{CnP-MIP}   & Q          & 0          & 27.22 (0.66)      & 0             & 1473.46     & 6.50          & 8.7           & 4.0              & 1.2            & 3.5              \\
			\textbf{CnP-MIP}   & Q          & 1          & 29.61 (0.61)      & 0             & 1476.85     & 4.20          & 14.0          & 2.0              & 0.5            & 11.5             \\
			\textbf{CnP-PATH}  & F          & -1         & 1.04              & 0             & 1327.47     & 12.50         & 19.2          & 6.3              & 12.9           & 0.0              \\
			\textbf{CnP-PATH}  & F          & 0          & 0.08              & 0             & 1325.23     & 6.40          & 8.1           & 3.4              & 1.6            & 3.1              \\
			\textbf{CnP-PATH}  & F          & 1          & 0.07              & 0             & 1361.74     & 4.60          & 15.0          & 2.2              & 0.5            & 12.3             \\
			\hline
			\texttt{n=2 m=40}  &            &            &                   &               &             &               &               &                  &                &                  \\
			\hline
			\textbf{SGM}       & -          & -          & 1.26              & 0             & 1348.56     & 13.70         & -             & -                & -              & -                \\
			\textbf{CnP-MIP}   & Q          & -1         & 27.87 (5.11)      & 0             & 1433.13     & 16.70         & 21.9          & 11.1             & 10.8           & 0.0              \\
			\textbf{CnP-MIP}   & Q          & 0          & 25.58 (3.53)      & 0             & 1434.09     & 12.80         & 13.4          & 8.2              & 1.1            & 4.1              \\
			\textbf{CnP-MIP}   & Q          & 1          & 29.72 (2.16)      & 0             & 1405.30     & 10.50         & 18.7          & 6.4              & 0.7            & 11.6             \\
			\textbf{CnP-PATH}  & F          & -1         & 0.89              & 0             & 1355.26     & 16.80         & 20.7          & 9.5              & 11.2           & 0.0              \\
			\textbf{CnP-PATH}  & F          & 0          & 0.70              & 0             & 1355.01     & 10.00         & 9.9           & 7.1              & 0.8            & 2.0              \\
			\textbf{CnP-PATH}  & F          & 1          & 0.62              & 0             & 1355.21     & 7.80          & 14.1          & 5.1              & 0.3            & 8.7              \\
			\bottomrule
		\end{tabular}
	}
\end{table}

\begin{table}[H]
	\centering
	\caption{\glspl{IPG} complete results, second set.  \label{tab:IPGRes2}}
	\resizebox{0.95\textwidth}{!}{
		\begin{tabular}{lll@{\hspace{4em}}rrrrrrrr}
			\toprule
			\textbf{Algorithm} & \textbf{O} & \textbf{C} & \textbf{Time (s)} & \textbf{\#TL} & \textbf{SW} & \textbf{\#It} & \textbf{Cuts} & \textbf{ESOCuts} & \textbf{VCuts} & \textbf{MIPCuts} \\
			\hline
			\texttt{n=3 m=40}  &            &            &                   &               &             &               &               &                  &                &                  \\
			\hline
			\textbf{SGM}       & -          & -          & 27.04             & 2             & 2339.79     & 20.10         & -             & -                & -              & -                \\
			\textbf{CnP-MIP}   & Q          & -1         & 140.33 (5.49)     & 0             & 2991.76     & 20.20         & 28.5          & 13.2             & 15.3           & 0.0              \\
			\textbf{CnP-MIP}   & Q          & 0          & 128.74 (3.06)     & 0             & 3016.22     & 11.60         & 15.6          & 8.9              & 1.9            & 4.8              \\
			\textbf{CnP-MIP}   & Q          & 1          & 162.20 (2.58)     & 0             & 2980.69     & 9.30          & 21.9          & 6.7              & 0.9            & 14.3             \\
			\textbf{CnP-PATH}  & F          & -1         & 2.35              & 0             & 2882.45     & 17.60         & 24.9          & 12.6             & 12.3           & 0.0              \\
			\textbf{CnP-PATH}  & F          & 0          & 0.87              & 0             & 2906.33     & 10.80         & 14.0          & 8.8              & 1.4            & 3.8              \\
			\textbf{CnP-PATH}  & F          & 1          & 0.79              & 0             & 2898.04     & 9.00          & 21.1          & 6.6              & 0.8            & 13.7             \\
			\hline
			\texttt{n=2 m=80}  &            &            &                   &               &             &               &               &                  &                &                  \\
			\hline
			\textbf{SGM}       & -          & -          & 14.97             & 1             & 2676.52     & 19.40         & -             & -                & -              & -                \\
			\textbf{CnP-MIP}   & Q          & -1         & 29.83 (11.47)     & 0             & 3127.96     & 7.60          & 6.7           & 5.4              & 1.3            & 0.0              \\
			\textbf{CnP-MIP}   & Q          & 0          & 27.02 (7.27)      & 0             & 3127.97     & 7.80          & 7.0           & 5.3              & 0.7            & 1.0              \\
			\textbf{CnP-MIP}   & Q          & 1          & 36.71 (10.06)     & 0             & 3124.63     & 6.10          & 8.6           & 3.6              & 0.5            & 4.5              \\
			\textbf{CnP-PATH}  & F          & -1         & 7.71              & 0             & 2914.36     & 8.80          & 8.1           & 6.7              & 1.4            & 0.0              \\
			\textbf{CnP-PATH}  & F          & 0          & 5.45              & 0             & 2926.82     & 7.00          & 6.1           & 4.5              & 0.4            & 1.2              \\
			\textbf{CnP-PATH}  & F          & 1          & 4.93              & 0             & 2936.52     & 5.80          & 7.4           & 3.4              & 0.4            & 3.6              \\
			\hline
			\texttt{n=2 m=100} &            &            &                   &               &             &               &               &                  &                &                  \\
			\hline
			\textbf{SGM}       & -          & -          & 77.13             & 3             & 2861.20     & 21.10         & -             & -                & -              & -                \\
			\textbf{CnP-MIP}   & Q          & -1         & 102.57 (36.29)    & 0             & 3750.38     & 10.30         & 10.9          & 7.4              & 3.5            & 0.0              \\
			\textbf{CnP-MIP}   & Q          & 0          & 105.97 (33.07)    & 1             & 3454.41     & 14.30         & 14.5          & 9.4              & 1.2            & 3.9              \\
			\textbf{CnP-MIP}   & Q          & 1          & 107.04 (30.86)    & 0             & 3771.62     & 12.00         & 18.0          & 6.3              & 0.8            & 10.9             \\
			\textbf{CnP-PATH}  & F          & -1         & 23.02             & 0             & 3496.86     & 11.22         & 11.67         & 8.33             & 3.33           & 0.0              \\
			\textbf{CnP-PATH}  & F          & 0          & 14.46             & 0             & 3488.44     & 10.70         & 11.0          & 7.1              & 1.2            & 2.7              \\
			\textbf{CnP-PATH}  & F          & 1          & 14.56             & 0             & 3507.71     & 10.30         & 14.8          & 6.4              & 0.7            & 7.7              \\
			\bottomrule
		\end{tabular}
	}
\end{table}

\end{document}